\documentclass[11pt, 3p]{elsarticle}

\usepackage{lineno,hyperref}
\usepackage{amstext}
\usepackage{amsmath,amssymb,amsthm,graphicx}
\modulolinenumbers[5]

\journal{Computers \& Operations Research}

\DeclareMathOperator*{\argmax}{arg\,max}
\DeclareMathOperator*{\argmin}{arg\,min}

\newtheorem{theorem}{Theorem}[section]

\newtheorem{lemma}[theorem]{Lemma}

\theoremstyle{remark}

\usepackage[short]{optidef}
\usepackage{changepage}
\usepackage[utf8x]{inputenc}
\usepackage[T1]{fontenc}
\usepackage{lmodern}
\usepackage[english]{babel}
\usepackage{hyperref}
\usepackage{comment}
\usepackage{enumitem}
\usepackage{algpseudocode}
\usepackage[ruled,vlined]{algorithm2e}
\usepackage{subcaption}
\usepackage{pgfplots}
\usepackage{booktabs}
\usepackage{tabularx}
\usepackage{ltablex}
\usepackage{longtable}
\usepackage{comment}
\usepackage{url}
\usepackage{bbm}
\pgfplotsset{compat=1.9}

\newcommand{\inprod}[2]{{\langle #1,#2 \rangle}} 

\newcommand{\Tcal}{\mathcal{T}}

\newcommand{\Scal}{\mathcal{S}}
\newcommand{\sn}{{\mathcal S}^n}

\newcommand{\R}{\mathbb{R}}

\newcommand{\textif}{\textrm{if}}

\newcommand{\trace}{\textrm{trace}}
\newcommand{\Diag}{\textrm{Diag}}

\biboptions{authoryear}

\usepackage{background}
\usepackage{setspace}
\begin{document}

\SetBgContents{Published at \url{https://doi.org/10.1016/j.cor.2022.105958}}      
\SetBgPosition{current page.center}
\SetBgAngle{0}                                    
\SetBgColor{gray}                                 
\SetBgScale{1.5}                                  
\SetBgHshift{0}                                   
\SetBgVshift{9cm} 

\begin{frontmatter}
\title{An Exact Algorithm for Semi-supervised Minimum Sum-of-Squares Clustering}


\author{Veronica Piccialli\corref{cor1}}
\ead{veronica.piccialli@uniroma1.it}
\address{Department of Computer, Control and Management Engineering, \\ Sapienza University of Rome, Via Ariosto 25, 00185, Italy}
\author{Anna Russo Russo}
\ead{anna.russo.russo@uniroma2.it}
\address{Department of Civil Engineering and Computer Science Engineering, \\ University of Rome Tor Vergata, Via del Politecnico 1, 00133, Italy}
\author{Antonio M. Sudoso}
\ead{antonio.maria.sudoso@uniroma2.it}
\address{Department of Civil Engineering and Computer Science Engineering, \\ University of Rome Tor Vergata, Via del Politecnico 1, 00133, Italy}
\cortext[cor1]{Corresponding author}

\begin{abstract}
The minimum sum-of-squares clustering (MSSC), or $k$-means type clustering, is traditionally considered an unsupervised learning task. In recent years, the use of background knowledge to improve the cluster quality and promote interpretability of the clustering process has become a hot research topic at the intersection of mathematical optimization and machine learning research. The problem of taking advantage of background information in data clustering
is called semi-supervised or constrained clustering. In this paper, we present branch-and-{cut} algorithm for semi-supervised MSSC, where background knowledge is incorporated as pairwise must-link and cannot-link constraints. For the lower bound procedure, we solve the semidefinite programming relaxation of the MSSC discrete optimization model, and we use a cutting-plane procedure for strengthening the bound. For the upper bound, instead, by using integer programming tools, we use an adaptation of the $k$-means algorithm to the constrained case.
For the first time, the proposed global optimization algorithm efficiently manages to solve real-world instances up to 800 data points with different combinations of must-link and cannot-link constraints and with a generic number of features. This problem size is about four times larger than the one of the instances solved by state-of-the-art exact algorithms.
\end{abstract}


\end{frontmatter}

\section{Introduction}
Cluster analysis or clustering is the task of grouping a set of patterns or observations, represented as points in a multidimensional space, in such a way that objects in the same group are more similar to each other than to those in other groups \citep{rao1971cluster, hartigan1979algorithm, jain1999data}. The groups are called clusters and the set of groups is a clustering.
Clustering is a fundamental tool in modern applications, such as pattern recognition, data mining, computer vision, machine learning, and knowledge discovery.
The clustering process is by definition unsupervised, which means that it only uses unlabeled data. However, without any supervision, clustering algorithms often produce solutions that are not relevant to expert knowledge. In order to improve the cluster quality and obtain meaningful solutions, researchers have focused on integrating knowledge to allow guidance on the clustering process. Incorporating such knowledge into the clustering process transforms the MSSC from an unsupervised learning problem to a semi-supervised one. Semi-supervised learning is a branch of machine learning where some information on the ground truth, as for example the labels of a subset of data, is available \citep{zhu2009introduction}.
 
Among many criteria used in cluster analysis, the most natural, intuitive, and frequently adopted criterion is the minimum sum-of-squares clustering (MSSC) or $k$-means type clustering. Given a set of $n$ data points in the $d$-dimensional Euclidean space, $p_1, \dots, p_n, \ p_i \in \mathbb{R}^d$, the MSSC aims to partition them into $k$ clusters, $C_1, \dots, C_k$, by minimizing the sum of squared Euclidean distances between the data points and the centers $m_1, \dots, m_k$ of the clusters they belong to. It can be formulated as:
\begin{equation}\label{eq:MSSC_unc}
\textrm{min} \ \sum_{j=1}^{k} \sum_{p_i \in C_j} \|p_i - m_j\|_2^2.
\end{equation}

In recent years, the MSSC has attracted a lot of attention in the area of data science and operations research \citep{gambella2021optimization}, and it has been extended to integrate background knowledge on objects and/or on clusters through user constraints \citep{basu2008constrained, brieden2017constrained, ganccarski2020constrained}. Nevertheless, with the presence of user constraints, clustering problems become harder and require the development of dedicated algorithms. The literature about solving the MSSC in the semi-supervised setting goes under the name of ``constrained clustering''. 

Constraints encountered in clustering problems can be categorized as instance-level and cluster-level constraints. Instance-level constraints are expressed as pairwise constraints, typically must-link and cannot-link, indicating that two points must, or cannot, be placed in the same cluster. These constraints can be inferred from class labels: if two objects have the same label then they are linked by a must-link constraint, otherwise by a cannot-link constraint. Pairwise constraints naturally arise in many domains such as gene clustering \citep{pensa2008constrained, maraziotis2012semi, tran2021local}, land consolidation \citep{borgwardt2014geometric}, and document clustering \citep{huang2006text, hu2008towards}. On the other hand, cluster-level constraints impose some knowledge on the structure of the clusters \citep{lai2021neighborhood, gnagi2021matheuristic}.
Constrained clustering methods can also be classified into two main categories: search-based (also known as constraint-based) methods and distance-based (also known as similarity-based) methods \citep{dinler2016survey}.
In search-based methods \citep{basu2004active, davidson2005clustering, zhang2019framework}, a clustering algorithm is modified to incorporate the prior knowledge into the clustering search space. 
This can be realized either by enforcing user constraints (hard clustering) or by including penalty terms for unsatisfied constraints in the objective function (soft clustering).
In distance-based methods \citep{xiang2008learning, xiaocui2020deep}, a clustering algorithm is used in conjunction with a novel distance measure modified in accordance with the prior knowledge expressed by the user constraints. For example, the distance measure should be adjusted in such a way that two observations involved in a must-link constraint will be closer to each other while two observations involved in a cannot-link constraint will be farther away from each other.
Furthermore, some hybrid methods integrate search-based and distance-based methods \citep{bilenko2004hybrid}.
In the following, we will focus on hard search-based methods, since our method falls into this class.

It is worth noticing that enforcing constraints in the context of an existing clustering algorithm is not a trivial step since, for a generic set of constraints, it is not straightforward to even determine whether a feasible clustering exists.
Feasibility analysis for instance-level and cluster-level constraints has received a lot of attention in the literature: \cite{davidson2005clustering, davidson2007intractability} discuss complexity results for the feasibility problem in presence of must-link and cannot-link constraints showing that determining the feasibility of general subsets of constraints is NP-complete.

Unconstrained MSSC is NP-hard \citep{pacheco2005scatter, aloise2009complexity}; as a consequence, exact algorithms for semi-supervised clustering are not very common in the literature. Moreover, these algorithms are incredibly sensitive to the number of observations and constraints involved, meaning that the majority of the existing methods are only capable to solve problems of very limited size. 

In this paper, we propose a global optimization algorithm for semi-supervised MSSC where prior knowledge on the data points is incorporated in the form of instance-level constraints. This method is based on the branch-and-{cut} technique, which is a very general framework for finding optimal solutions in discrete optimization problems. Lower bounds are obtained by relaxing the discrete optimization model for semi-supervised MSSC with pairwise constraints and solving the resulting Semidefinite Programming (SDP) relaxation. 
The main contributions of the paper are:
\begin{enumerate}
    \item {We solve semi-supervised MSSC to global optimality with an SDP-based branch-and-cut algorithm}.
    \item We propose an SDP-based heuristic for solving semi-supervised MSSC providing a bound on the optimality gap of the produced clustering.
    \item Numerical experiments show that our semi-supervised algorithm solves instances up to $n=800$ points with up to $n/2$ cannot-link and/or must-link constraints, that is a problem size about 4 times larger than the one solved by state-of-the-art exact algorithms. 
\end{enumerate}

\section{Related Work}
This section presents a review of search-based clustering methods with pairwise constraints, both heuristic and exact. 
As for the heuristic approaches, most of them have the following properties in common: (1) they extend an unconstrained clustering algorithm to integrate user constraints; (2) they integrate instance-level constraints (i.e.{,} must-link and cannot-link constraints) (3) they are usually fast and find an approximate solution, and therefore do not guarantee the satisfaction of all the constraints nor the global solution even when the constraints are satisfied. 

The most notable example of a search-based method enforcing pairwise hard constraints on the MSSC problem is COP-$k$-means \citep{wagstaff2001constrained}. COP-$k$-means adapts the classic $k$-means algorithm to handle must-link and cannot-link constraints: in each iteration it tries to assign each observation to the nearest cluster center so that no pairwise constraints are violated, but performs this task in a greedy fashion without contemplating any backtracking.
As a consequence, not only this algorithm lacks optimality guarantees but it can also fail to return a solution when a feasible assignment exists. Indeed, clustering solutions produced by COP-$k$-means depend both on the initial assignment, as in the classic $k$-means algorithm, and on the assignment order in each iteration.
To address the issue of constraint violation in COP-$k$-means, \cite{tan2010icop} and \cite{rutayisire2011sequenced} propose a modified version comprising a pre-computation of the assignment order: in the former algorithm, named ICOP-$k$-means, the assignment order is based on the certainty of each data point, computed through the technique of clustering ensemble; in the latter, the order is determined by carrying out a breadth-first search of the cannot-link set.
\cite{huang2008semi} propose MLC-KMeans and use the concept of assistant centroids, which are calculated using the points implicated by must-link constraints for each cluster. The assistant centers are also used to compute the similarity of data points and clusters.
Two recent heuristics make use of integer programming. In \cite{vrain2020constrained}, integer programming is used to post-process the solution produced by an unconstrained algorithm and force the pairwise constraints. \cite{baumann2020binary} defines a $k$-means like heuristic, where the assignment step solves a binary programming problem to enforce all the pairwise constraints.

Some clustering algorithms relax the user constraints, and consequently do not guarantee to satisfy all of them. Methods falling in this class use penalties as a trade-off between finding the best cluster assignment and satisfying as many constraints as possible.
\cite{basu2004active} propose PCKmeans, a clustering method with pairwise constraints that allows some constraints to be violated. PCKmeans minimizes a modified version of the MSSC objective function by adding penalties for the violation of must-link and cannot-link constraints yielding a soft-clustering assignment.
\cite{davidson2005clustering} propose another soft constrained clustering algorithm. This method penalizes constraint violation using the Euclidean distance: if a must-link constraint is violated, then the cost is given by the distance between the centroids of the clusters containing the data points that should be together; if a cannot-link constraint is violated, then the penalty is computed as the distance between the center of the cluster the two data points are assigned to and the nearest cluster center.
The Lagrangian constrained clustering approach in \cite{ganji2016lagrangian} considers a penalty for violating only cannot-link constraints. Must-link constraints are used to aggregate the data points into super-points so that they are all satisfied. This method uses an iterative approach where in each iteration a Lagrangian relaxation is solved by increasing penalties for constraints which remain unsatisfied in subsequent iterations. Lastly, \cite{gonzalez2020dils} propose an iterative local search procedure for clustering instances with a big number of pairwise constraints.

\medskip

Although there has been considerable methodological research activity in the area of semi-supervised clustering, there exist few exact methods for semi-supervised MSSC and most of them are extensions of unconstrained MSSC exact algorithms. All these approaches use general optimization tools, such as integer programming or constraint programming, and they search for a global optimum that satisfies all the constraints but can only solve instances with limited number of data points.
Global optimization algorithms proposed in the literature for unconstrained MSSC are based on cutting plane \citep{sherali2005global, xia2005cutting}, branch-and-bound \citep{koontz1975branch, diehr1985evaluation, brusco2006repetitive, krislock2016computational}, branch-and-cut \citep{aloise2009branch, piccialli2021sos} and column generation algorithms \citep{du1999interior, aloise2012improved}.

To the best of our knowledge, the first exact method for semi-supervised MSSC is presented in \cite{xia2009global}. Xia extends the global optimization method in \cite{xia2005cutting} for unsupervised MSSC to deal with instance-level constraints. Approximate results are obtained by halting the algorithm before global convergence, but according to \cite{aloise2012improved}, this kind of method can produce exact solutions only on instances with about 25 data points.
\cite{aloise2012column} extend their previous column-generation algorithm \citep{aloise2012improved} for unsupervised MSSC and show that they are able to solve instances with less than 200 data points and few tens of must-link and cannot-link constraints.
The column generation algorithm in \cite{aloise2012column} is also used as a starting point in \cite{babaki2014constrained}, where the authors propose a column-generation framework to solve the semi-supervised MSSC with must-link and cannot-link constraints, as well as other monotonic constraints. However, its application scope is limited to instances with less than 200 data points.
Most of the recently proposed approaches are based on the constraint programming paradigm which is a general framework for solving combinatorial optimization problems; among these, \cite{duong2013declarative, duong2015constrained, duong2017constrained} describe a declarative framework for several optimization criteria, including that of MSSC.
\cite{duong2015constrained} show that their method outperforms the column-generation algorithm in \cite{babaki2014constrained}, but results are reported only on instances consisting of less than 200 data points. 
\cite{guns2016repetitive} extend the repetitive branch-and-bound algorithm (RBBA) proposed by Brusco for unconstrained MSSC \citep{brusco2006repetitive}. They replace the internal branch-and-bound of RBBA by a constraint programming solver, and use it to compute tight lower and upper bounds. This algorithm, named as CPRBBA, can incorporate must-link and cannot-link constraints and can be used in a multi-objective constrained clustering setting, which minimizes the MSSC objective and maximizes the separation between clusters. However, results are presented only for few instances of less than 200 data points with up to 250 constraints.

The hardest task when developing exact algorithms for constrained MSSC is to compute good bounds in a reasonable amount of time. 
Although their focus is not on exact methods, \cite{liberti2021side} propose several MINLP reformulations of MSSC with side constraints, both exact and approximate. Among these, some reformulations yield convex continuous relaxations that may be embedded in global optimization algorithms.

In the next section we describe the ingredients of the proposed branch-and-{cut} algorithm. 

\section{Branch-and-Cut Algorithm}
\subsection{Notation}
Let $\mathcal{N} = \{1, \dots, n\}$ be the set of indices of the data points $p_i \in \mathbb{R}^d$, $\mathcal{K} = \{1, \dots, k\}$ the set of indices of the clusters, $\mathcal{ML} \subseteq \mathcal{N} \times \mathcal{N}$ the set of must-link (ML) constraints and $\mathcal{CL} \subseteq \mathcal{N} \times \mathcal{N}$ the set of cannot-link (CL) constraints.
Let $\sn$ denote the set of all $n\times n$ real symmetric matrices. We denote by $M\succeq 0$ that matrix $M$ is positive semidefinite and let ${\mathcal S}_+^n$ be the set of all positive semidefinite matrices of size
$n\times n$. We denote by $\inprod{\cdot}{\cdot}$ the 
trace inner product. That is, for any
$A, B \in \mathbb{R}^{m\times n}$, we define $\inprod{A}{B}:= \trace (B^\top A )$. Given a matrix $A$, we denote by $A_{i\cdot}$ the $i$-th row of $A$. We denote by $e_n$ the vector of all ones of length $n$ and by $I_n$ the identity matrix of size $n \times n$. We omit the subscript in case the dimension is clear from the context. For a symmetric matrix $A$ we denote the set of its eigenvalues by $\lambda(A)$.

\subsection{Problem definition}
The semi-supervised MSSC with pairwise constraints can be formulated as follows:
\begin{mini!}[2]
{}{\sum_{i=1}^{n} \sum_{j=1}^{k} x_{ij} \| p_i - m_j \|_2^2}
    {\label{prob:semiMSSC}}{}
\addConstraint{\sum_{j=1}^k x_{ij}}{= 1 \label{con:hard}}{\quad \forall i \in \mathcal{N}}
\addConstraint{\sum_{i=1}^{n} x_{ij}}{\geq \label{con:cardinality} 1}{\quad \forall j \in \mathcal{K}}
\addConstraint{x_{ih}}{= x_{jh} \label{con:ml}}{\quad \forall h \in \mathcal{K}, \ \forall(i,j) \in \mathcal{ML}}
\addConstraint{x_{ih} + x_{jh}}{\leq 1 \label{con:cl}}{\quad \forall h \in \mathcal{K}, \ \forall (i,j) \in \mathcal{CL}}
\addConstraint{x_{ij}}{\in \{0, 1\}}{\quad \forall i \in \mathcal{N}, \ \forall j \in \mathcal{K}}
\addConstraint{m_j}{\in \mathbb{R}^d}{\quad \forall j \in \mathcal{K}}.
\end{mini!}

In \eqref{prob:semiMSSC}, the cluster centers are at the unknown points $m_j$ and the binary decision variable $x_{ij}$ expresses whether data point $i$ is assigned to cluster $j$ or not. Constraints (\ref{con:hard}) ensure that each data point is assigned to exactly one cluster, Constraints (\ref{con:cardinality}) avoid empty clusters and Constraints (\ref{con:ml}) and (\ref{con:cl}) enforce ML and CL constraints, respectively. From now on, we denote the discrete feasible region of \eqref{prob:semiMSSC} by $\mathcal{F} = \big\{x_{ij} \in \{0, 1\} \ \big\vert \ \eqref{con:hard}, \eqref{con:cardinality}, \eqref{con:ml}, \eqref{con:cl}  \ \textrm{hold} \big\}$.


Problem \eqref{prob:semiMSSC} is a mixed-integer programming problem with nonlinear objective function, which is known to be NP-hard like the unconstrained MSSC \citep{davidson2007intractability, aloise2009complexity}. The difficulty of this problem derives from two different aspects. First, the decision variables are discrete and the pairwise constraints are embedded as hard constraints. Secondly, the objective is nonlinear and nonconvex. These difficulties in the objective as well as in the constraints make semi-supervised MSSC extremely hard and challenging to solve exactly.

Setting the gradient of the objective function with respect to $m_j$ to zero yields
\begin{equation*}
m_j = \frac{\sum_{i=1}^n x_{ij} p_i}{\sum_{i=1}^n x_{ij}}, \quad \forall j \in \mathcal{K},
\end{equation*}
that is the average of all the points assigned to cluster $j$. In other words, the optimal cluster centers are always at the centroids of the clusters.
Replacing the formula for $m_j$ in~\eqref{prob:semiMSSC}, we obtain 
\begin{mini}[2]
{}{\sum_{i=1}^n \sum_{j=1}^k x_{ij}\Bigg\|p_i - \frac{\sum_{l=1}^n x_{lj} p_l}{\sum_{l=1}^n x_{lj}}\Bigg\|_2^2}
    {\label{prob:semiMSSC2}}{}
\addConstraint{x_{ij}}{\in \mathcal{F}}{\quad \forall i \in \mathcal{N}, \ \forall j \in \mathcal{K}}.
\end{mini}
The unconstrained MSSC, i.e.{,} Problem \eqref{prob:semiMSSC2} with  $\mathcal{ML}=\mathcal{CL}=\emptyset$, is known to be equivalent to a nonlinear SDP problem, the so-called 0-1 SDP \citep{peng2007approximating}. Here we derive first the equivalent 0-1 SDP in the unconstrained setting, following the derivation process in \cite{peng2007approximating}, and then we extend the formulation to the constrained case. 

For convenience, we collect all the data points $p_i$ as rows in a matrix $W_p$. Let $W = W_p W_p^\top$ be the matrix of the inner products of the data points, i.e.{,} $W_{ij} = p_i^\top p_j$ for all $i,j \in \mathcal{N}$. Furthermore, collect the binary decision variables $x_{ij}$ from \eqref{prob:semiMSSC2} in the $n\times k$ matrix $X$ and define the $n \times n$ matrix $Z := X(X^\top X)^{-1}X^\top.$ 
It is easy to verify that $Z$ is a symmetric matrix with nonnegative entries satisfying $Z^2 = Z$. Furthermore, it follows immediately that the sum of each row (or column) is equal to 1 and the trace is equal to $k$. By rearranging the terms in the objective function we obtain the equivalent reformulation of unconstrained MSSC:
\begin{mini}[2]
{}{\textrm{tr}(W(I - Z))}
    {\label{prob:PengSDP}}{}
\addConstraint{Ze=e, \ \textrm{tr}(Z) = k}{}{}
\addConstraint{Z \ge 0, \ Z^2 = Z, \ Z=Z^\top}{}{}.
\end{mini}

\cite{peng2007approximating} proved that any feasible solution $Z$ for Problem \eqref{prob:PengSDP} is
necessarily associated {with} a feasible unconstrained MSSC assignment matrix $X$. Therefore, the equivalence among the MSSC formulations \eqref{prob:semiMSSC} (without constraints \eqref{con:ml} and \eqref{con:cl}) and \eqref{prob:PengSDP} is established. Note that for a given data matrix $W_p$, the trace of $W$ is a fixed quantity. Therefore, we can solve Problem \eqref{prob:PengSDP} by maximizing $\textrm{tr}(WZ)$.
From now on we refer to any feasible solution $Z$ as ``clustering matrix''. For any feasible clustering matrix it is easy to verify that:

\begin{itemize}
    \item If data points $i$ and $j$ are in the same cluster C, then $Z_{i\cdot} = Z_{j\cdot}$ and the non-zero entries in rows $i$ and $j$ are equal to $\frac{1}{|C|}$, where $|C|$ is the cardinality of cluster C.
    \item If data points $i$ and $j$ are not in the same cluster, then $Z_{ij} = 0$.
\end{itemize}

This allows to express ML and CL constraints by equality constraints on the matrix $Z$, getting the following 0-1 SDP reformulation of semi-supervised MSSC with pairwise constraints:
\begin{mini}[2]
{}{\textrm{tr}(W(I - Z))}
    {\label{prob:SemiPengSDP}}{}
\addConstraint{Ze=e, \ \textrm{tr}(Z) = k}{}{}
\addConstraint{Z_{ih}}{= Z_{jh}}{\quad \forall h \in \mathcal{N}, \ \forall(i,j) \in \mathcal{ML}}
\addConstraint{Z_{ij}}{= 0}{\quad \forall(i,j) \in \mathcal{CL}}
\addConstraint{Z \ge 0, \ Z^2 = Z, \ Z=Z^\top}{}{}.
\end{mini}

It is worthwhile comparing \eqref{prob:semiMSSC2} with \eqref{prob:SemiPengSDP}. Differently from \eqref{prob:semiMSSC2}, the objective function in \eqref{prob:SemiPengSDP} is
linear and the constraints are all linear except $Z^2 = Z$, that is even more complex than the binary constraints in \eqref{prob:semiMSSC2} since it constrains the eigenvalues of $Z$ being binary. However, by relaxing this constraint, we can get very good lower bounds on the original problem, as shown in the next subsection. The SDP relaxation we get is the main ingredient of our branch-and-cut approach.

\subsection{Lower bound}
By relaxing the non-convex constraint $Z^2=Z$ in \eqref{prob:SemiPengSDP} into the semidefinite constraint $Z \succeq 0$,  we obtain the following SDP relaxation:
\begin{mini}[2]
{}{\textrm{tr}(W(I - Z))}
    {\label{prob:SemiPengSDPRelax}}{}
\addConstraint{Ze=e, \ \textrm{tr}(Z)= k}{}{}
\addConstraint{Z_{ih}}{= Z_{jh}}{\quad \forall h \in \mathcal{N}, \ \forall(i,j) \in \mathcal{ML}}
\addConstraint{Z_{ij}}{= 0}{\quad \forall(i,j) \in \mathcal{CL}}
\addConstraint{Z}{\ge 0, \ Z \in \mathcal{S}^n_+}{}.
\end{mini}

Interestingly enough, ML constraints allow to reduce the size of the SDP relaxation to be solved at each node. In particular, we construct a graph from the set of ML constraints where each node is a data point, and there is an edge between two nodes if the corresponding data points are involved in a must-link constraint. Then we compute the transitive closure of the graph leading to a partition of nodes into $s$ components $B_1, \dots, B_s$. 
These connected components collect data points that must be in the same cluster. We collapse each connected component $B_i$ into a single ``super point'' $\overline{p}_i$ that is obtained as the sum of the data points in $B_i$. Therefore, instead of clustering the set of initial points $p_1, \dots, p_n$, we search for a clustering on the set of super points $\overline{p}_1, \dots, \overline{p}_s$. A similar approach for handling ML constraints, and thus reducing the size of the problem, is used in \cite{guns2016repetitive} where they replace all the points in each connected component with a single point which is the sum of all the points in that set, weighted by the size of the connected component. In our case, we exploit this idea to obtain an equivalent reformulation over a lower dimensional positive semidefinite cone.

Formally, given a set of $n$ initial data points, assume that there are $s \leq n$ connected components $B_1, \dots, B_s$ with super points $\overline{p}_1, \dots, \overline{p}_s$.
Define the $s \times n$ binary matrix
\[
\Tcal^{s}_{ij} = \left\{ \begin{array}{ll}
1 & \textif~ j \in B_i\\
0 & \textif~ j \notin B_i
\end{array}\right.
\]
encoding the data points that have been merged (i.e., the $i$-th row of $\mathcal{T}^s$ is the indicator vector of $B_i$ for $i=1,\dots,s$)
and vector $e^s = \Tcal^{s} e$ containing the size of each connected component. Since we are interested in clustering the super points arising from ML constraints, we observe that $\Tcal^{s} W (\Tcal^{s})^\top$ shrinks the size of $W_p$ (i.e.{,} the matrix having data points as rows) by substituting row $i$ by $\overline{p}_i = \sum_{p \in B_i} p$ and omitting the rows of the merged data points. 
Instance-level constraints that remain to be satisfied are only CL constraints. To this end, we denote by $\overline{\mathcal{CL}}$ the set of CL constraints between two connected components. We add to $\overline{\mathcal{CL}}$ a CL between super-points $\overline{p}_i$ and $\overline{p}_j$ if there exists a CL constraint on two data points $p$ and $q$ such that $p \in B_i$ and $q \in B_j$. Therefore, the SDP solved for computing the lower bound becomes
\begin{mini}[2]
{}{\inprod{I}{W} - \inprod{\Tcal^{s} W (\Tcal^{s})^\top}{Z}}
    {\label{prob:SDPbab}}{}
\addConstraint{Z e^{s}}{= e }{}
\addConstraint{\inprod{\Diag(e^{s})}{Z}}{= k}{}
\addConstraint{Z_{ij}}{= 0}{\quad \forall(i,j) \in \overline{\mathcal{CL}}}
\addConstraint{Z \ge 0, \ Z \in \Scal_+^{s}}{}{}.
\end{mini}

The proof of equivalence between \eqref{prob:SemiPengSDPRelax} and \eqref{prob:SDPbab} can be easily derived by using Theorem 2 in \cite{piccialli2021sos}. In the latter, ML and CL constraints are added one at a time when visiting the branch-and-bound tree, since the children are generated either by merging two points thanks to {an} ML or adding the corresponding CL constraint. Therefore, the size of the matrix is reduced by one at each branching level, and the effect of the size reduction is essentially negligible thanks to the small number of visited nodes.
In the semi-supervised setting, instead, pairwise constraints are included at the beginning of the formulation, leading to a different transformation (the one described by $\Tcal^{s}$) and to a significant reduction in size. In particular, for the case $k=2$ we can further reduce the size of SDP since a CL constraint determines in which cluster a data point should be. Hence, we can derive additional ML constraints as follows: for any $i, j, h \in \mathcal{N}$ such that $(i, h) \in \mathcal{CL}$ and $(j, h) \in \mathcal{CL}$, then we add $(i, j)$ to $\mathcal{ML}$.

Our algorithm produces a binary enumeration tree. Every time a node is split into two children, a pair of points $(i, j)$ is chosen, and {an} ML constraint and a CL constraint are respectively imposed on the left and the right child. These two types of constraints partition the set of solutions associated to
the parent child into two disjoint subsets. As for the branching strategy, we observe that in a matrix $Z$ corresponding to a clustering, for each pair $(i,j)$ either $Z_{ij}=0$ or $Z_{i\cdot} = Z_{j\cdot}$. Suppose that for the optimal solution of the SDP relaxation there are indices $i$ and $j$ such that $Z_{ij} \| Z_{i\cdot}-Z_{j\cdot} \|_2^2 \neq 0$, then we generate a cannot-link branch with $Z_{ij} = 0$ and a must-link branch with $Z_{i\cdot} = Z_{j\cdot}$. Regarding the variable selection, we choose indices $i$ and $j$ such that in both branches we expect a significant improvement  of the lower bound. The branching pair is chosen as the \[\argmax_{i, j} \Big\{ \min \big\{Z_{ij}, \| Z_{i\cdot}-Z_{j\cdot} \|_2^2 \big\} \Big\}.\]

\subsection{Valid inequalities}
In this subsection, we present three classes of inequalities we use for strengthening the bound.

\paragraph{Pair} In any feasible clustering, it holds that
\begin{equation}\label{eq:pairs}
        Z_{ij}\le Z_{ii},\quad Z_{ij}\le Z_{jj}\quad \forall i,j \in \mathcal{N}, \ i\not=j.
\end{equation} 

\paragraph{Triangle} If data points $i$ and $j$ are in the same cluster and data points $j$ and $h$ are in the same cluster, then $i$ and $h$ necessarily must be in the same cluster:
\begin{equation}\label{eq:triangle}
            Z_{ij}+Z_{ih}\le Z_{ii}+Z_{jh}\quad \forall i,j,h \in \mathcal{N}, \ i,j,h ~\mathrm{distinct}.
\end{equation}

\paragraph{Clique} If the number of clusters is $k$, for any subset $Q$ of $k+1$ points at least two points must be in the same cluster:
\begin{equation}\label{eq:clique}
        \sum_{(i,j)\in Q,i<j}Z_{ij}\ge \frac{1}{n-k+1} \quad\forall Q\subset\mathcal{N},\ |Q|=k+1.
\end{equation}
  
\cite{piccialli2021sos} report detailed computational results on benchmark instances showing that lower bounds provided by the Peng-Wei SDP relaxation with these inequalities are very close to the optimal values. Here we add pair, triangle and clique inequalities to the SDP relaxation via a cutting-plane procedure only if they are violated. 
After each cutting-plane iteration we remove constraints that are not active at the optimal solution of the SDP relaxation. In this way, we keep each bounding problem to a computationally tractable size.
Moreover, inequalities that are included in the parent model during the last cutting-plane iteration are also inherited by its children and added in their model from the beginning. This procedure allows to quickly retrieve several effective inequalities, and save a significant number of cutting-plane iterations in the children.

\subsection{Valid lower bounds and feasibility issues}
Off-the-shelf solvers for SDP problems are generally based on Interior Point Methods (IPMs) \citep{alizadeh1995interior}. 
While these algorithms solve SDPs to arbitrary accuracy in polynomial time, in practice, they suffer from scalability and performance issues. When the considered clustering problems have large number of variables and constraints, solution time and memory requirements tend to explode, making SDPs prohibitive to work with. Compared to IPMs, solvers based on semi-proximal ADMM/ALM (Alternating Direction Method of Multipliers/Augmented Lagrangian Method) can scale to significantly larger problem sizes, while trading off the accuracy of the resulting output \citep{sun2015convergent, yang2015sdpnal}. 
However, they show two limitations that need to be addressed to guarantee the theoretical validity of a branch-and-bound algorithm.
First, when using ADMM/ALM methods, it is hard to reach a solution to high precision in a reasonable amount of time. This implies that the bound provided by an SDP relaxation solved to a moderate precision may not be a valid lower bound.
Second, existing SDP solvers based on ADMM/ALM frameworks can not detect primal and dual infeasibilities since they are designed on the assumption that the primal and the dual problems are feasible. 

\paragraph{Valid lower bounds}
To overcome the first limitation, we use the post-processing technique developed in \cite{JaChayKeil2007} where the authors propose a method to obtain rigorous lower bounds on the 
optimal value of SDP problems. Recall that the dual objective function value of any dual feasible solution yields a bound on the optimal objective function of the primal. Therefore, every dual feasible solution and, in particular, the optimal solution of the dual SDP, gives a valid bound on the solution of the discrete optimization problem. However, the dual objective function value represents a valid dual bound only if the SDP relaxation is solved to high precision. In the following, we generate a ``safe'' underestimate for the primal SDP, and hence a valid lower bound for the MSSC problem with pairwise constraints. We start with the following lemma.
\begin{lemma}\label{lemma:jansson}
Let $S, X \in \mathcal{S}^n$ be matrices that satisfy
$0 \leq \lambda_{\min}(X)$ and $\lambda_{\max}(X) \leq \bar{x}$ for some $\bar{x} \in \R$.
Then the inequality
\begin{equation*}
    \inprod{S}{X} \geq \bar{x}\sum_{i \colon  \lambda_i(S) <0}\lambda_i(S)
\end{equation*}
holds.
\end{lemma}
\begin{proof}
See the proof in \cite[Lemma 3.1 on p. 184]{JaChayKeil2007}.
\end{proof}
At this point, we can adapt the post-processing procedure in \cite{JaChayKeil2007} to our SDP. Let $E_i^s$ be the symmetric matrix such that $\inprod{E_i^s}{Z} = (Ze^s)_i$ for $i \in \{1, \dots, s\}$ and $E_{ij}$ be symmetric matrix such that $\inprod{E_{ij}}{Z} = Z_{ij}$ for $(i,j) \in \overline{\mathcal{CL}}$. Furthermore, {we} define a mapping $\tau$ such that for $(i, j) \in \overline{\mathcal{CL}}$ we have $\tau(i, j) = h \in \{1, \dots, |\mathcal{\overline{CL}}|\}$.
\begin{theorem}\label{theorem:pp}
Let $Z^\star$ be the optimal solution of Problem \eqref{prob:SDPbab} with objective function value $p^\star$. Consider the dual variables $\tilde{y} \in \mathbb{R}^{s+1}$, $\tilde{\mu} \in \mathbb{R}^{|\mathcal{\overline{CL}}|}$, $\tilde{V} \in \mathcal{S}^s$, $\tilde{V} \geq 0$ and set $\tilde{S} = -\mathcal{T}^s W (\mathcal{T}^s)^\top - \sum_{i=1}^{s} \tilde{y}_i E_i^s - \tilde{y}_{s+1} \Diag(e^s) - \sum_{(i, j) \in \mathcal{\overline{CL}}} \tilde{\mu}_{\tau(i,j)} E_{ij}  - \tilde{V}$, then a safe lower bound for $p^\star$ is given by
\begin{equation*}
    lb = \sum_{i=1}^{s} \tilde{y}_i + k \tilde{y}_{s+1} + \bar{z} \sum_{i\colon \lambda_i(\tilde{S}) < 0} \lambda_i(\tilde{S}),
\end{equation*}
where $\bar{z} \geq \lambda_{\max}(Z^\star)$.
\end{theorem}
\begin{proof}
See the proof in \ref{app:proof}.
\end{proof}

If matrix $\tilde{S}$ is positive semidefinite, then $(\tilde{y}, \tilde{\mu}, \tilde{V}, \tilde{S})$ is a dual feasible solution and $\sum_{i=1}^{s} \tilde{y}_i + k \tilde{y}_{s+1}$ is already a valid lower bound. Otherwise, we decrease the dual objective by adding the negative perturbation $\bar{z} \sum_{i\colon \lambda_i(\tilde{S}) < 0} \lambda_i(\tilde{S})$ to it where $\bar{z}$ is an upper bound on the maximum eigenvalue of any feasible solution of the primal. We can bound the maximum eigenvalue of any feasible $Z$ by 1 since $Z e \leq Z e^s = e$. Therefore, we set $\bar{z} = 1$ and after the SDP relaxation has been solved to a moderate precision, we apply Theorem \ref{theorem:pp} to obtain valid lower bounds.

\paragraph{Feasibility}
It is well known that for infeasible convex optimization problems some of the iterates of ADMM/ALM methods diverge. However, terminating the algorithm when the iterates become large is unreliable in practice: such termination criterion is just an indication that a problem might be infeasible, and not a certificate of infeasibility.
Given a large number of initial pairwise constraints and those derived from branching decisions, it is very likely to visit infeasible nodes while going down the tree. Hence, we want to detect infeasibility as soon as possible to save ourselves the expensive part of calculating the bound.
To this end, by using an off-the-shelf integer programming solver {such as Gurobi \citep{gurobi}}, we solve the feasibility Problem (\ref{prob:feasibility}) to check whether a feasible clustering exists. If such clustering does not exist then we simply prune the node, otherwise we solve the SDP relaxation to obtain a lower bound.
\begin{mini}[2]
{}{0}{\label{prob:feasibility}}{}
\addConstraint{x_{ij}}{\in \mathcal{F}}{\quad \forall i \in \mathcal{N}, \ \forall j \in \mathcal{K}}.
\end{mini}

\subsection{Heuristic}
\label{section:heuristic}
The most popular heuristic for solving unconstrained MSSC is unarguably the $k$-means algorithm \citep{macqueen1967some, lloyd1982least}. Given the initial cluster centers, $k$-means proceeds by alternating between two steps until convergence: in the first step, each data point is assigned to the closest cluster center, whereas in the second step, the cluster centers are updated by taking the average of all the data points that belong to each cluster. These two steps are repeated until the assignment of data points to clusters no longer changes.


Because of its simplicity, efficiency and empirical success, it has been commonly used as a template for developing constrained clustering algorithms. In this direction, COP-$k$-means adapts $k$-means to handle ML and CL constraints: in each iteration it tries to assign each data point to the nearest cluster center so that no pairwise constraints are violated \citep{wagstaff2001constrained}.
\cite{vrain2020constrained} propose a novel approach for constrained clustering developing a post-processing procedure to enforce constraints a posteriori. Given a matrix that contains the degree of membership of each data point to each cluster (obtained by means of any constrained or unconstrained clustering algorithm) this method finds a hard assignment satisfying all the constraints by solving a combinatorial optimization problem. The post-processed clustering matrix does not minimize the within-cluster sum of squares criterion (i.e.{,} the MSSC objective function) but it is constructed in such a way that it resembles the initial cluster membership. In \cite{vrain2020constrained}, clustering solutions provided by this procedure score a better accuracy than those obtained by running COP-$k$-means several times, though no optimality guarantee is given.
Another recent heuristic algorithm for constrained clustering has been proposed by \cite{baumann2020binary}. The main idea of this approach consists in modifying the assignment step of $k$-means to find an optimal clustering satisfying all the pairwise constraints by using an Integer Linear Programming (ILP) formulation. Differently from $k$-means, this assignment step does not admit an analytic solution, but it requires the solution of a linear program with binary variables. The objective function and the feasible set are the same as in \eqref{prob:semiMSSC}, but the centers are fixed at the current value.  After the assignment step, this procedure behaves exactly like $k$-means, since the cluster centers are updated by averaging the data points assigned to each of them. The initial cluster centers are randomly picked from the set of data points.
Unlike COP-$k$-means, if there exists a feasible clustering the assignment step can never fail to assign each data point to a cluster.

In order to make our branch-and-bound algorithm efficient, we need a heuristic procedure to obtain a feasible high quality clustering at each node. Greedy heuristic algorithms tend to fail when the number of CL and ML is high, while $k$-means, also in its constrained variant proposed by \cite{baumann2020binary}, is very sensitive to the choice of the initial cluster centers, similarly to other local methods for non-convex optimization problems. To overcome this drawback, several initialization strategies have been proposed to prevent $k$-means to get stuck in a low quality local minimum \citep{pena1999empirical, celebi2013comparative}.
In \cite{piccialli2021sos}, an initialization technique has been defined based on the primal solution provided by the SDP relaxation. The idea is that if the SDP relaxation \eqref{prob:SDPbab} were tight, then the solution $Z$ would be feasible for the 0-1 SDP \eqref{prob:SemiPengSDP}, and hence would allow to easily recover the centroids. If the relaxation is not tight, then we find the closest rank-$k$ approximation of $Z$ in terms of the Frobenius norm. More {precisely}, we solve an optimization problem where the norm of the difference between $Z$ and a rank-$k$ approximating matrix is minimized. This problem admits an analytic solution that can be obtained by computing the truncated singular value decomposition (SVD) of $Z$: instead of taking all the eigenvalues and their corresponding eigenvectors, we only take the $k$ largest eigenvalues and their corresponding eigenvectors (see the low-rank approximation theorem in \cite{eckart1936approximation}).

In our branch-and-cut algorithm we use an effective heuristic combining this SDP-based initialization technique with the assignment step developed by \cite{baumann2020binary}. We name this procedure IPC-$k$-means, {short for integer pairwise constrained $k$-means}, and we describe it in Algorithm \ref{alg:ipckmeans}.

\begin{algorithm}
\caption{IPC-$k$-means}
\label{alg:ipckmeans}

\KwIn{Data points $p_1, \dots, p_n$, number of clusters $k$, sets of constraints $\mathcal{ML}$ and $\mathcal{CL}$, data matrix $W_p$, optimal solution $\tilde{Z}$ of the SDP relaxation with $\mathcal{ML}$ and $\mathcal{CL}$ constraints}

\begin{enumerate}[leftmargin=*, nolistsep]
\item Solve $\hat{Z} = \argmin \{||\tilde{Z} - Z||_F \ \textrm{s.t.} \ \textrm{rank}(Z) =  k\}$ by computing the truncated SVD of $\tilde{Z}$.
\item Compute the centroid matrix approximation $\hat{M} = \hat{Z}W_p$.
\item Cluster the rows of $\hat{M}$ with $k$-means to get the initial cluster centers $m_1, \dots, m_k$.
\end{enumerate}
\Repeat{convergence}{
\begin{enumerate}[leftmargin=*, rightmargin=30pt, nolistsep]
    \item[4.1.] Compute the optimal cluster assignments $x_{ij}^\star$ by solving:
    \begin{mini}[2]
    {}{\sum_{i=1}^{n} \sum_{j=1}^{k} x_{ij} \| p_i - m_j \|_2^2 \label{prob:ass}}
        {}{}
    \addConstraint{x_{ij}}{\in \mathcal{F}}{\quad \forall i \in \mathcal{N}, \ \forall j \in \mathcal{K}}.
    \end{mini}
    \item[4.2.] Set $C_j \leftarrow \{ p_i : x^\star_{ij} = 1\}$ for each $j=1,\dots,k$.
    \item[4.3.] Update the cluster centers $m_1, \dots, m_k$ by taking the mean of the data points assigned to each cluster $C_1, \dots, C_k$.
\end{enumerate}
}
\KwOut {Clusters $C_1, \dots, C_k$.}

\end{algorithm}

Our heuristic requires the solution of a small number of SDP problems. Hence it is more expensive from the computational point of view in comparison to random initialization techniques. However, solving the SDP relaxation provides a lower bound on the optimal value, so that it allows to certify an optimality gap for the provided feasible clustering. Furthermore, numerical experiments in Section \ref{sec:numresheur} show that the initial choice of cluster centers retrieved from the SDP relaxation that incorporates ML and CL constraints is robust and yields high quality clustering solutions.
The overall branch-and-cut algorithm is illustrated in Algorithm \ref{alg:bbpseudocode}. 

\begin{algorithm}
\caption{Branch-and-Cut Algorithm}
\label{alg:bbpseudocode}

\KwIn{Sets of constraints $\mathcal{ML}$ and $\mathcal{CL}$, number of clusters $k$, linear kernel matrix $W$}

\begin{enumerate}[label*=\arabic*., leftmargin=*, nolistsep] 
    \item Build $\mathcal{T}^s$, $e^s$ and $\mathcal{\overline{CL}}$ from $\mathcal{ML}$ and $\mathcal{CL}$. Let $P_0$ be the initial 0-1 SDP in (\ref{prob:SemiPengSDP}) and set $\mathcal{Q} = \{P_0\}$.
    \item Set $X^\star = \textrm{null}$ with objective value $v^\star = \infty$.
    \item While $\mathcal{Q}$ is not empty:
    \begin{enumerate}[label*=\arabic*., rightmargin=30pt, nolistsep]
        \item Select and remove problem $P$ from $\mathcal{Q}$.
        \item Solve the feasibility problem (\ref{prob:feasibility}). If it is infeasibile, go to \textit{Step 3}.
        \item Solve the SDP relaxation (\ref{prob:SDPbab}) to get a lower bound $LB$ and the optimal solution $Z$.
        \item If $LB \geq v^\star$, go to \textit{Step 3}.
        \item Search for pair (\ref{eq:pairs}), triangle (\ref{eq:triangle}), and clique (\ref{eq:clique}) inequalities violated by $Z$. If any are found, add them to the current SDP relaxation and go to \textit{Step 3.3}.
        \item Run the heuristic in Algorithm \ref{alg:ipckmeans} to get an assignment $X$ and an upper bound $UB$. If $UB < v^\star$ then set $v^\star \leftarrow UB$, $X^\star \leftarrow X$.
        \item Select the branching pair $(i,j)$ and partition problem $P$ into ML and CL sub-problems. For each problem update $\mathcal{T}^s$, $e^s$ and $\mathcal{\overline{CL}}$ accordingly, add them to $\mathcal{Q}$ and go to \textit{Step 3}.
    \end{enumerate}
\end{enumerate}

\KwOut{Optimal assignment matrix $X^\star$ with objective value $v^\star$}

\end{algorithm}


\section{Computational Results}
In this section, we describe the details of the implementation and we show the numerical results on real-world datasets.

\subsection{Implementation details}
PC-SOS-SDP, which stands for Pairwise Constrained SOS-SDP, is implemented in C++ with some routines written in MATLAB. The SDP relaxation at each node is solved by means of SDPNAL+, a MATLAB software that implements an augmented Lagrangian method to solve large scale SDPs with bound constraints \citep{sun2020sdpnal+}. We set the accuracy tolerance of the solver to $10^{-5}$ in the relative KKT residual. We also use Gurobi \citep{gurobi} in order to deal with the ILPs required for the upper bound computation and the feasibility check. We run the experiments on a machine with Intel(R) Xeon(R) 8124M CPU @ 3.00GHz with 16 cores, 64 GB of RAM, and Ubuntu Server 20.04. For pair and triangle inequalities, we randomly separate at most $100000$ valid cuts, we sort them in decreasing order with respect to the violation, and we add the first 5\% of violated ones in the current cutting-plane iteration. For the separation of clique inequalities, we use the heuristic procedure described in \cite{piccialli2021sos} that returns at most $n$ valid cuts. 
The tolerance for checking the violation is set to $10^{-4}$. The maximum number of cutting plane iterations at the root node is set to 50, whereas for the ML and CL children this number is set to 30. We stop the cutting-plane procedure not only when there are no violated inequalities, but also when the lower bound does not improve significantly in the current cutting-plane iteration. Finally, we visit the tree with the best-first search strategy.
In order to improve the efficiency of the branch-and-bound search, PC-SOS-SDP processes many nodes in parallel using a thread pool of fixed size: whenever an ML or CL sub-problem is created, it is assigned to one of the available threads and run in parallel with the other threads of the pool. Furthermore, each thread calls SDPNAL+ in a separate MATLAB session. 
For the parallel setting, we use different configurations depending on the instance size. For small instances ($n < 300$) we use a pool of 16 threads, each of them running on a MATLAB session with a single component thread. For larger instances ($n \geq 300$) we use a pool of 8 threads, each of them running on a MATLAB session with 2 component threads. In all cases, the session for the computation at the root node uses all the available cores. The source code of PC-SOS-SDP and the instances used in our tests are available at \url{https://github.com/antoniosudoso/pc-sos-sdp}.

\subsection{Instances generation}
We build our semi-supervised clustering instances from real-world datasets for classification problems.
For each dataset, we generate several instances differing in the type and amount of user constraints. Specifically, given a dataset with $n$ data points, we chose to build sets of constraints consisting of approximately $n/2$ and $n/4$ constraints; each set either contains only ML constraints, only CL constraints, or an equal number of ML and CL constraints. All these constraints are enforced according to the dataset true class partitioning, which is, in general, in contradiction with the unconstrained MSSC global optimum.
For each of these configurations, we generate 5 random sets of constraints using a classic procedure described in \cite{wagstaff2001constrained} and then also used in more recent works on exact semi-supervised clustering methods \citep{babaki2014constrained, duong2015constrained, guns2016repetitive}: at each step, a pair of data points is randomly selected and either {an} ML or a CL constraint is defined depending on the true labels of the data points; the procedure stops when the desired amount of ML and CL constraints is achieved. 

\subsection{Results on real-world instances}
We consider 12 real-world datasets for classification problems, with a number of data points $n$ ranging between $150$ and $801$, and with a number of features $d$ ranging between $4$ and $20531$. 
For each dataset, we consider 10 instances with only CL constraints, 10 instances with only ML constraints, and 10 instances with an equal amount of ML and CL constraints. Overall, we build and solve 360 constrained clustering instances. 
In the experiments, an instance is solved successfully when the optimality gap is less or equal than $10^{-4}$. This gap measures the difference between the best upper and lower bounds and it is calculated as $(UB - LB) / UB$.
The MSSC requires the user to specify the number of clusters $k$ to generate. In the literature, clustering validity indices are commonly used to determine a suitable number of clusters. In the semi-supervised setting instead, the number of clusters is known and assumed to be equal to the number of classes. The datasets characteristics, i.e.{,} number of data points, features, and clusters, are reported in Table~\ref{tab:datasets}.

\begin{table}
\centering
\begin{tabular}{lccc}
\toprule
Dataset   &  $n$   &  $d$ & $k$  \\
\midrule
 Iris & 150 & 4 & 3 \\ 
 Wine & 178 & 13 & 3  \\
 Connectionist & 208 & 60 & 2 \\
 Seeds & 210 & 7 & 3 \\
 Glass & 214 & 9 & 6 \\
 Heart & 299 & 12 & 2 \\
\bottomrule
\end{tabular}\hspace{20pt}
\begin{tabular}{lccc}
\toprule
Dataset   &  $n$   &  $d$ & $k$  \\
\midrule
Vertebral & 310 & 6 & 2 \\
Accent & 329 & 12 & 6 \\
Ecoli & 336 & 7 & 8 \\
ECG5000 & 500 & 140 & 5 \\
Computers & 500 & 720 & 2 \\
Gene & 801 & 20531 & 5\\
\bottomrule
\end{tabular}
\caption{Characteristics of the real-world datasets. They all can be downloaded at the UCI \citep{uci} and UCR \citep{UCRArchive2018} websites.}
\label{tab:datasets}
\end{table}

The results of our experiments are outlined in Tables \ref{tab:ml}, \ref{tab:cl}, \ref{tab:mlcl}, comprising respectively tests on instances with only ML constraints, only CL constraints and both types of constraints. Every entry of these tables involves a single dataset, whose name, size and number of clusters are reported, and shows aggregated statistics of 5 random instances with a certain number of ML and CL constraints: these statistics include the average number of separate data points at the root node (i.e.{,} the size of the problem at the root), the average root gap, the average size of the branching tree and the average completion time. Finally, in the last column, a percentage expresses the share of successfully solved instances, i.e.{,} instances that were solved before reaching the maximum number of nodes, which is set to 200; in those cases where one or more instances can not be solved within this limit, in the last column, between brackets, we report the average gap reached before halting the branch-and-{cut} algorithm on the unsolved instances.

\begin{figure}
    \centering
    \includegraphics[scale = 0.75]{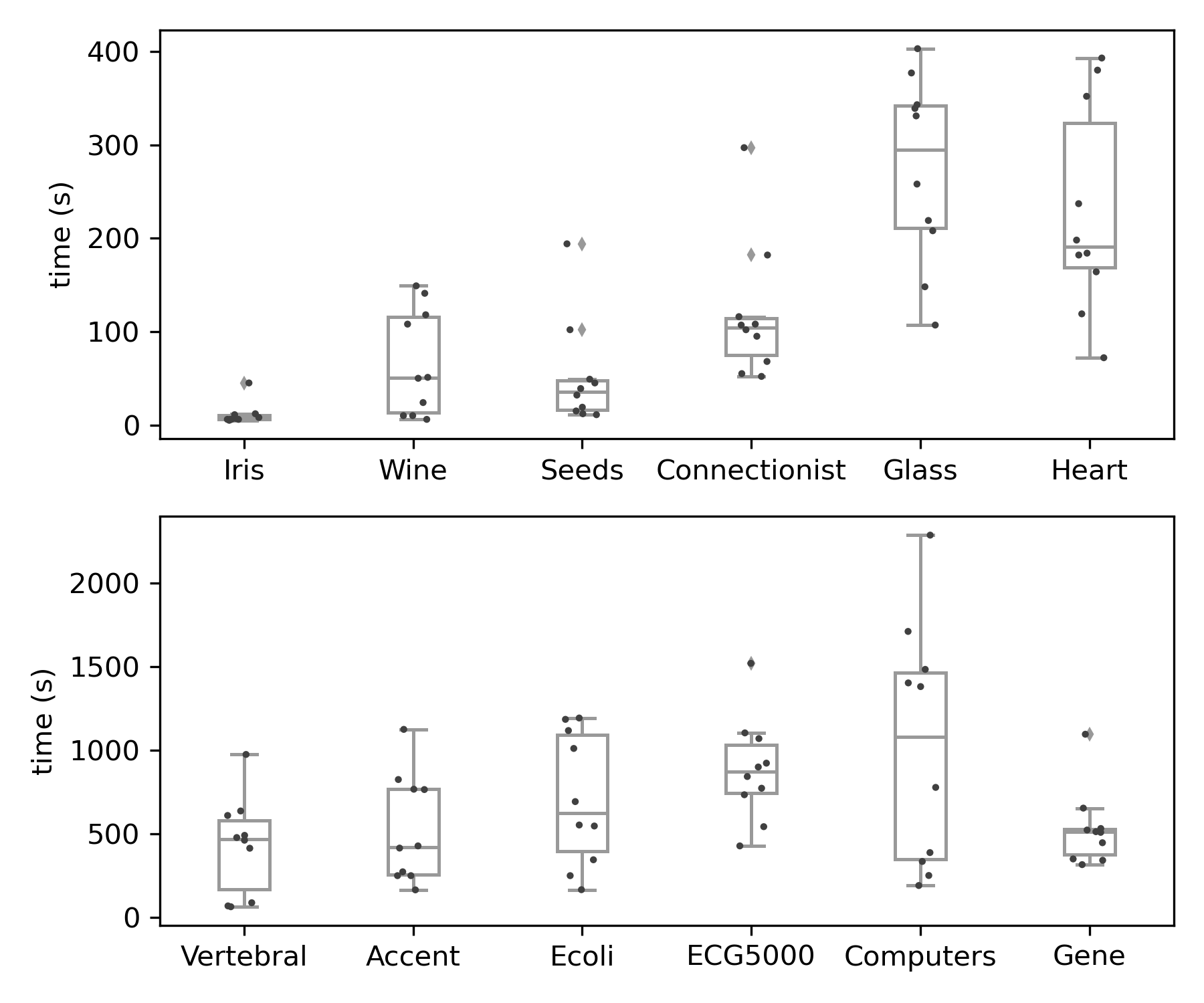}
    \caption{Boxplot of computational times on instances with ML constraints only.}
    \label{fig:tempi_ml}
\end{figure}

\paragraph{Must-Link constraints only} Table \ref{tab:ml} shows that, when dealing with ML constraints only, our algorithm can solve to optimality every instance in less than half an hour on average and within a handful of nodes. It can be noted that the mean number of nodes and mean completion time are generally lower when a higher number of constraints are included: indeed, adding {an} ML constraint has the effect of merging two separate data points into one, thus overall decreasing the initial size of the clustering problem. The boxplots in Figure \ref{fig:tempi_ml} offer a more detailed view of the computational time required on each of our 120 instances with ML constraints only, grouped by dataset. Here we can see that on datasets with up to 300 points our branch-and-{cut} algorithm always converged to optimality in less than 500 seconds, while on bigger instances the highest time required is slightly over 2000 seconds, with $95\%$ of these instances being solved in less than 1500 seconds.

\begin{figure}
    \centering
    \includegraphics[scale = 0.75]{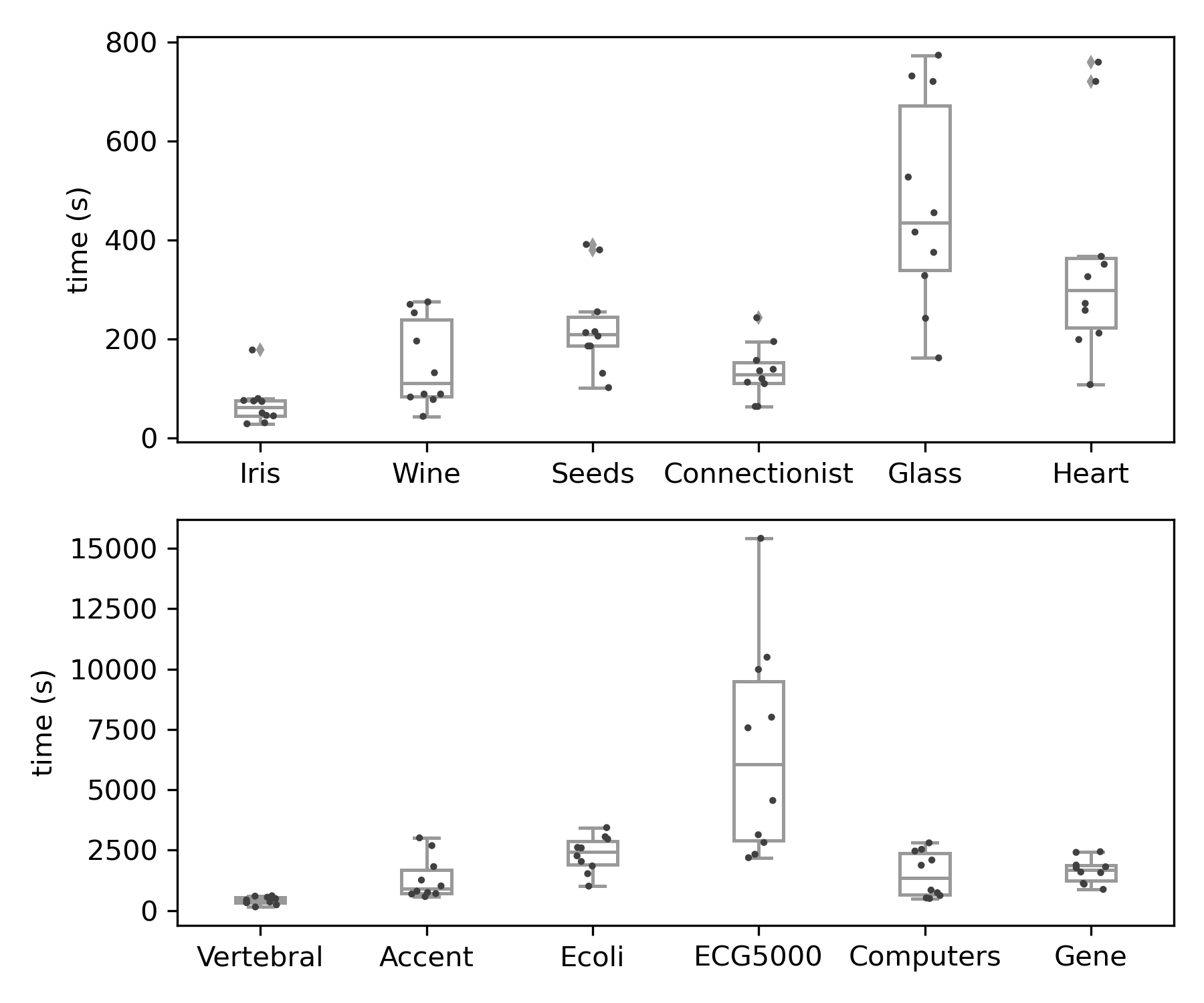}
    \caption{Boxplot of computational times on instances with CL constraints only.}
    \label{fig:tempi_cl}
\end{figure}

\paragraph{Cannot-Link constraints only} Table \ref{tab:cl} displays results for instances consisting of CL constraints only. These kind{s} of constraints usually make the clustering problem much harder than the unconstrained version. Indeed, the computational time required to solve these instances to optimality tends to grow larger as a greater number of constraints is included. Instances with just two clusters represent a consistent exception to this tendency since in this case CL constraints allow to infer non-redundant ML constraint, therefore decreasing the size of the initial clustering problem. Overall, 4 out of 120 instances can not be solved within the threshold that is set on the branching tree size, but for these instances we are still able to provide a very good clustering solution with a certified relative gap smaller than $0.03\%$. As can be seen in Figure \ref{fig:tempi_cl}, the instances with less than 300 points are all solved in 800 seconds or less; for the bigger instances the maximum computational time spent is about 16000 seconds but $91\%$ of them are solved to optimality in less than 2 hours.

\begin{figure}
    \centering
    \includegraphics[scale = 0.75]{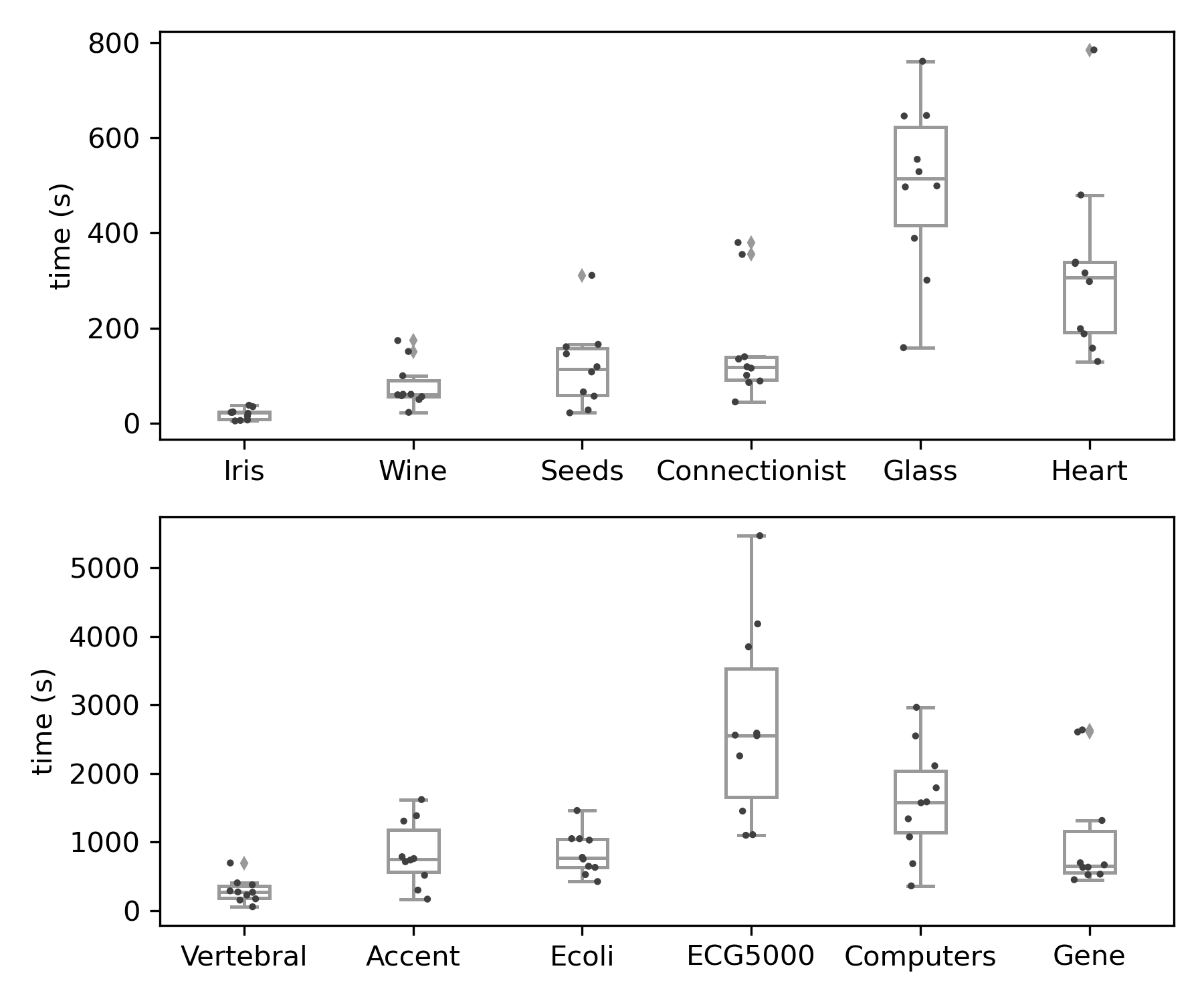}
    \caption{Boxplot of computational times on instances with ML and CL constraints.}
    \label{fig:tempi_mlcl}
\end{figure}

\paragraph{Combination of Must-Link and Cannot-Link constraints} Results for instances with mixed types of constraints are reported in Table \ref{tab:mlcl}. When dealing with both types of constraints, our algorithm fails to solve 2 out of 120 instances, as the maximum number of nodes is encountered before reaching optimality. Nonetheless on these instances we find a feasible solution with relative gap not greater than $0.04\%$. All the remaining instances are solved within one hour on average: the highest time required was 800 seconds for datasets made of less than 300 data points, and barely more than 5000 seconds for larger datasets, as shown in the boxplots of Figure \ref{fig:tempi_mlcl}. 

Note that the time needed to solve the ILPs for finding feasible clusterings is neglectable: on the largest instance, it is {about 10 seconds}.

These results show that our method is able to solve successfully instances up to a size of $n = 801$ data points and $n/2$ pairwise constraints, with CL constraints being the most challenging kind. Moreover, we can not miss to point out that our average root gap is smaller than $1\%$ on each dataset and for each type of constraints, and smaller than $0.01\%$ on $47\%$ of our instances.
Furthermore, it is worth noticing that we are able to solve instances with a very large number of features $d$ (over 20 thousand), as our algorithm is minimally sensitive to the dimension of the feature space. Indeed, the number of features has no influence on the lower bound computation since it is hidden in the matrix $W$ that is computed only once.

Non-aggregated statistics on each of our instances can be found in \ref{app:num}, where each PC-SOS-SDP execution is described in more detail.

\begin{table}
\centering
\begin{tabular}{ccccccccc}
\toprule
            dataset &  constraints &   $n$ &  $k$ &  size & root gap &  nodes &  time (s) & success rate \\
\midrule
               Iris &           50 & 150 &  3 & 100.4 &   0.004\% &    1.8 &      16.6 &         100\% \\
               Iris &          100 & 150 &  3 &  56.0 &   0.001\% &    1.0 &       5.8 &         100\% \\
               Wine &           50 & 178 &  3 & 128.2 &   0.016\% &    1.8 &     101.8 &         100\% \\
               Wine &          100 & 178 &  3 &  81.4 &   0.009\% &    1.4 &      31.6 &         100\% \\
              Seeds &           50 & 210 &  3 & 160.0 &   0.032\% &    4.6 &      84.4 &         100\% \\
              Seeds &          100 & 210 &  3 & 111.6 &   0.009\% &    1.8 &      19.2 &         100\% \\
      Connectionist &           50 & 208 &  2 & 158.2 &   0.005\% &    1.0 &     160.8 &         100\% \\
      Connectionist &          100 & 208 &  2 & 109.0 &   0.007\% &    1.8 &      75.6 &         100\% \\
              Glass &           50 & 214 &  6 & 164.4 &   0.038\% &    7.0 &     234.2 &         100\% \\
              Glass &          100 & 214 &  6 & 117.8 &   0.059\% &    4.6 &     312.4 &         100\% \\
              Heart &          100 & 299 &  2 & 199.0 &   0.001\% &    1.0 &     309.2 &         100\% \\
              Heart &          150 & 299 &  2 & 152.8 &   0.001\% &    1.0 &     147.0 &         100\% \\
          Vertebral &          100 & 310 &  2 & 210.6 &   0.018\% &    5.0 &     594.8 &         100\% \\
          Vertebral &          150 & 310 &  2 & 161.4 &   0.050\% &    4.6 &     260.8 &         100\% \\
             Accent &          100 & 329 &  6 & 229.4 &   0.049\% &    3.0 &     371.0 &         100\% \\
             Accent &          150 & 329 &  6 & 186.0 &   0.374\% &    5.0 &     679.6 &         100\% \\
              Ecoli &          100 & 336 &  8 & 236.4 &   0.059\% &   33.4 &     852.0 &         100\% \\
              Ecoli &          150 & 336 &  8 & 189.4 &   0.015\% &    3.4 &     558.2 &         100\% \\
            ECG5000 &          150 & 500 &  5 & 350.4 &   0.014\% &    9.0 &    1086.4 &         100\% \\
            ECG5000 &          250 & 500 &  5 & 253.2 &   0.324\% &    6.2 &     679.2 &         100\% \\
          Computers &          150 & 500 &  2 & 350.4 &   0.003\% &    1.4 &    1652.2 &         100\% \\
          Computers &          250 & 500 &  2 & 252.0 &   0.001\% &    1.0 &     387.6 &         100\% \\
               Gene &          200 & 801 &  5 & 601.6 &   0.001\% &    1.0 &     503.8 &         100\% \\
               Gene &          400 & 801 &  5 & 417.2 &   0.001\% &    1.0 &     550.6 &         100\% \\
\bottomrule
\end{tabular}

\caption{PC-SOS-SDP results on the instances where only ML constraints are included.}
\label{tab:ml}
\end{table}

\begin{table}
\centering
\begin{tabular}{ccccccccc}
\toprule
            dataset &  constraints &   $n$ &  $k$ &  size & root gap &  nodes &  time (s) & success rate (gap) \\
\midrule
               Iris &           50 & 150 &  3 & 150.0 &   0.010\% &    2.2 &      78.4 &         100\% \\
               Iris &          100 & 150 &  3 & 150.0 &   0.102\% &    2.6 &      58.6 &         100\% \\
               Wine &           50 & 178 &  3 & 178.0 &   0.039\% &    1.4 &      85.2 &         100\% \\
               Wine &          100 & 178 &  3 & 178.0 &   0.063\% &    6.2 &     216.6 &         100\% \\
              Seeds &           50 & 210 &  3 & 210.0 &   0.054\% &    9.0 &     227.2 &         100\% \\
              Seeds &          100 & 210 &  3 & 210.0 &   0.097\% &   11.0 &     225.8 &         100\% \\
      Connectionist &           50 & 208 &  2 & 186.6 &   0.002\% &    1.0 &     165.6 &         100\% \\
      Connectionist &          100 & 208 &  2 & 139.4 &   0.006\% &    1.4 &     102.6 &         100\% \\
              Glass &           50 & 214 &  6 & 214.0 &   0.104\% &   32.6 &     539.4 &         100\% \\
              Glass &          100 & 214 &  6 & 214.0 &   0.108\% &   23.4 &     406.4 &         100\% \\
              Heart &          100 & 299 &  2 & 239.6 &   0.001\% &    1.0 &     481.8 &         100\% \\
              Heart &          150 & 299 &  2 & 182.4 &   0.001\% &    1.0 &     232.6 &         100\% \\
          Vertebral &          100 & 310 &  2 & 254.8 &   0.008\% &    1.4 &     425.0 &         100\% \\
          Vertebral &          150 & 310 &  2 & 199.2 &   0.030\% &   12.2 &     387.0 &         100\% \\
             Accent &          100 & 329 &  6 & 329.0 &   0.067\% &   57.2 &    1264.4 &          80\% (0.013\%)\\
             Accent &          150 & 329 &  6 & 329.0 &   0.107\% &   75.0 &    1393.0 &         100\% \\
              Ecoli &          100 & 336 &  8 & 336.0 &   0.104\% &  150.6 &    2432.0 &          60\% (0.024\%)\\
              Ecoli &          150 & 336 &  8 & 336.0 &   0.270\% &  116.6 &    2231.8 &         100\% \\
            ECG5000 &          150 & 500 &  5 & 500.0 &   0.032\% &   64.8 &    8042.8 &          80\% (0.012\%)\\
            ECG5000 &          250 & 500 &  5 & 500.0 &   0.032\% &   22.6 &    5254.2 &         100\% \\
          Computers &          150 & 500 &  2 & 427.2 &   0.001\% &    1.0 &    2351.2 &         100\% \\
          Computers &          250 & 500 &  2 & 321.2 &   0.001\% &    1.0 &     642.4 &         100\% \\
               Gene &          200 & 801 &  5 & 801.0 &   0.011\% &    2.6 &    1540.0 &         100\% \\
               Gene &          400 & 801 &  5 & 801.0 &   0.014\% &    3.0 &    1768.0 &         100\% \\
\bottomrule
\end{tabular}

\caption{{PC-SOS-SDP results on the instances where only CL constraints are included.}}
\label{tab:cl}
\end{table}

\begin{table}
\centering
\begin{tabular}{ccccccccc}
\toprule
            dataset &  constraints &   $n$ &  $k$ &  size & root gap &  nodes &  time (s) & success rate (gap) \\
\midrule
               Iris &           25 & 150 &  3 & 125.0 &   0.014\% &    1.8 &      28.2 &         100\% \\
               Iris &           50 & 150 &  3 & 100.0 &   0.004\% &    1.4 &      11.2 &         100\% \\
               Wine &           25 & 178 &  3 & 153.0 &   0.016\% &    1.4 &      65.8 &         100\% \\
               Wine &           50 & 178 &  3 & 128.4 &   0.022\% &    2.2 &      93.0 &         100\% \\
              Seeds &           25 & 210 &  3 & 185.0 &   0.034\% &    5.8 &     164.6 &         100\% \\
              Seeds &           50 & 210 &  3 & 160.0 &   0.042\% &    3.8 &      72.2 &         100\% \\
      Connectionist &           25 & 208 &  2 & 177.6 &   0.005\% &    1.0 &     225.2 &         100\% \\
      Connectionist &           50 & 208 &  2 & 139.6 &   0.003\% &    1.0 &      88.0 &         100\% \\
              Glass &           25 & 214 &  6 & 189.0 &   0.168\% &   52.2 &     542.0 &         100\% \\
              Glass &           50 & 214 &  6 & 164.6 &   0.431\% &   53.4 &     454.6 &         100\% \\
              Heart &           50 & 299 &  2 & 231.2 &   0.021\% &    1.8 &     427.8 &         100\% \\
              Heart &           75 & 299 &  2 & 189.4 &   0.001\% &    1.0 &     218.0 &         100\% \\
          Vertebral &           50 & 310 &  2 & 244.0 &   0.010\% &    3.4 &     385.6 &         100\% \\
          Vertebral &           75 & 310 &  2 & 199.2 &   0.017\% &    3.8 &     200.8 &         100\% \\
             Accent &           50 & 329 &  6 & 279.2 &   0.050\% &   16.2 &     783.4 &         100\% \\
             Accent &           75 & 329 &  6 & 254.4 &   0.442\% &  103.4 &     878.4 &          60\% (0.03\%)\\
              Ecoli &           50 & 336 &  8 & 286.0 &   0.050\% &   27.0 &     760.8 &         100\% \\
              Ecoli &           75 & 336 &  8 & 261.4 &   0.295\% &   47.4 &     912.2 &         100\% \\
            ECG5000 &           75 & 500 &  5 & 425.2 &   0.044\% &   33.8 &    3338.6 &         100\% \\
            ECG5000 &          125 & 500 &  5 & 375.0 &   0.069\% &   29.0 &    2087.8 &         100\% \\
          Computers &           75 & 500 &  2 & 404.0 &   0.012\% &    1.4 &    2045.2 &         100\% \\
          Computers &          125 & 500 &  2 & 323.2 &   0.001\% &    1.0 &    1166.0 &         100\% \\
               Gene &          100 & 801 &  5 & 701.0 &   0.008\% &    2.6 &    1185.6 &         100\% \\
               Gene &          200 & 801 &  5 & 602.2 &   0.006\% &    2.6 &     957.4 &         100\% \\
\bottomrule
\end{tabular}

\caption{{PC-SOS-SDP results on the instances where both ML and CL constraints are included.}}
\label{tab:mlcl}
\end{table}

\bigskip
We directly compare the results of our algorithm to those of CPRBBA, developed by \cite{guns2016repetitive}.
CPRBBA is chosen as a benchmark because, to the best of our knowledge, it is the most recent and most efficient global algorithm for solving MSSC with pairwise constraints.
CPRBBA is run on the same instances and on the same machine used for testing our own algorithm. A time limit of one hour is set on each instance.
Tables \ref{tab:competitor_ml}, \ref{tab:competitor_cl} and \ref{tab:competitor_mlcl} show the results obtained by CPRBBA on the 4 smallest datasets that we considered in our experiments, i.e.{,} Iris, Wine, Connectionist and Seeds.
We also made an attempt to solve instances consisting of just ML constraints for the datasets Glass, Heart, Vertebral and Ecoli but none of those instances were solved within our time limit. \cite{guns2016repetitive} only report results for datasets with about 200 data points or less, therefore we feel safe in assuming that the CPRBBA approach would fail to find the optimal solution on all the instances that are missing in Tables \ref{tab:competitor_ml}, \ref{tab:competitor_cl} and \ref{tab:competitor_mlcl}.
In the same tables we compare the results obtained by CPRBBA with the performance of PC-SOS-SDP. More specifically, for each dataset and instance type we report the average resolution time and the success rate of both algorithms.
For CPRBBA we consider success as solving to optimality an instance within one hour of time; for our algorithm PC-SOS-SDP we define success as solving the instance to optimality within one hour of time and within 200 nodes.
As can be seen in these tables, despite being impressively efficient at solving instances of the Iris dataset when ML constraints are involved, CPRBBA does not behave well on slightly bigger instances: not only it fails to solve any of the Connectionist instances, but its success rate drastically falls whenever CL constraints are enforced.
On the other hand, in the same time frame our algorithm is able to solve all the considered instances to optimality. It also appears consistently faster than CPRBBA on almost all the instances, with the only exceptions being instances of the Iris dataset.
Moreover, even in those cases marked as a failure, our branch-and-{cut} approach is usually able to provide a very good clustering solution (maybe even the optimal one), whose quality is guaranteed by the extremely low gap obtained before halting the branch-and-cut procedure. On the contrary, CPRBBA can not produce a feasible clustering solution before the algorithm termination, since it works by incrementally adding data points to the constrained MSSC until the entire dataset is considered. 
We also stress that the average gap found by PC-SOS-SDP at the root node is already extremely low, and this excellent performance is due to the effectiveness of our heuristic procedure, as we can see in the next section, combined with the goodness of the lower bound.

\begin{table}
\centering
\begin{tabular}{cccccccc}
\toprule
      & & & &   \multicolumn{2}{c}{CPRBBA} & \multicolumn{2}{c}{PC-SOS-SDP} \\
            dataset &  constraints &   $n$ &  $k$ & time (s) & success rate & time (s) & success rate \\
\midrule
               Iris &           50 & 150 &  3 &     0.20 &         100\% &      16.6 &         100\% \\
               Iris &          100 & 150 &  3 &     0.02 &         100\% &       5.8 &         100\% \\
               Wine &           50 & 178 &  3 &   275.58 &         100\% &     101.8 &         100\% \\
               Wine &          100 & 178 &  3 &    39.02 &         100\% &      31.6 &         100\% \\
      Connectionist &           50 & 208 &  2 &  3600.00 &           0\% &     160.8 &         100\% \\
      Connectionist &          100 & 208 &  2 &  3600.00 &           0\% &      75.6 &         100\% \\
              Seeds &           50 & 210 &  3 &   174.08 &         100\% &      84.4 &         100\% \\
              Seeds &          100 & 210 &  3 &    35.76 &         100\% &      19.2 &         100\% \\
\bottomrule
\end{tabular}
\caption{Comparison between CPRBBA and PC-SOS-SDP on instances with ML constraints only.}
\label{tab:competitor_ml}
\vspace{20pt}
\centering
\begin{tabular}{cccccccc}
\toprule
      & & & &   \multicolumn{2}{c}{CPRBBA} & \multicolumn{2}{c}{PC-SOS-SDP} \\
            dataset &  constraints &   $n$ &  $k$ & time (s) & success rate & time (s) & success rate \\
\midrule
               Iris &           50 & 150 &  3 &   750.74 &          80\% &      78.4 &         100\% \\
               Iris &          100 & 150 &  3 &   724.54 &          80\% &      58.6 &         100\% \\
               Wine &           50 & 178 &  3 &  3600.00 &           0\% &      85.2 &         100\% \\
               Wine &          100 & 178 &  3 &  3600.00 &           0\% &     216.6 &         100\% \\
      Connectionist &           50 & 208 &  2 &  3600.00 &           0\% &     165.6 &         100\% \\
      Connectionist &          100 & 208 &  2 &  3600.00 &           0\% &     102.6 &         100\% \\
              Seeds &           50 & 210 &  3 &  3008.29 &          20\% &     227.2 &         100\% \\
              Seeds &          100 & 210 &  3 &  3600.00 &           0\% &     225.8 &         100\% \\
\bottomrule
\end{tabular}

\caption{Comparison between CPRBBA and PC-SOS-SDP on instances with CL constraints only.}
\label{tab:competitor_cl}
\vspace{20pt}
\centering
\begin{tabular}{cccccccc}
\toprule
      & & & &   \multicolumn{2}{c}{CPRBBA} & \multicolumn{2}{c}{PC-SOS-SDP} \\
            dataset &  constraints &   $n$ &  $k$ & time (s) & success rate & time (s) & success rate \\
\midrule
               Iris &           25 & 150 &  3 &     0.30 &         100\% &      28.2 &         100\% \\
               Iris &           50 & 150 &  3 &     3.32 &         100\% &      11.2 &         100\% \\
               Wine &           25 & 178 &  3 &  1642.63 &          60\% &      65.8 &         100\% \\
               Wine &           50 & 178 &  3 &  3136.97 &          20\% &      93.0 &         100\% \\
      Connectionist &           25 & 208 &  2 &  3600.00 &           0\% &     225.2 &         100\% \\
      Connectionist &           50 & 208 &  2 &  3600.00 &           0\% &      88.0 &         100\% \\
              Seeds &           25 & 210 &  3 &   782.40 &          80\% &     164.6 &         100\% \\
              Seeds &           50 & 210 &  3 &   458.63 &         100\% &      72.2 &         100\% \\
\bottomrule
\end{tabular}
\caption{Comparison between CPRBBA and PC-SOS-SDP on instances with both ML and CL constraints.}
\label{tab:competitor_mlcl}
\end{table}

\subsection{Heuristic algorithm results}\label{sec:numresheur}

We test the efficiency of our heuristic algorithm IPC-$k$-means by comparing the quality of its clustering solutions with the solutions produced by three other heuristic methods for semi-supervised clustering: COP-$k$-means, the post-processing approach (PP) proposed by \cite{vrain2020constrained} applied to the unconstrained clustering obtained by $k$-means, and the heuristic $BLPKM_{CC}$ proposed by \cite{baumann2020binary}.
These heuristic algorithms address the MSSC problem with ML and CL constraints while ensuring that no user constraints are violated in the final clustering partition.

All the competitor heuristic methods are executed with 100 different initializations, all built with the $k$-means++ algorithm  \citep{vassilvitskii2006k}, whereas IPC-$k$-means is run only once since it comprises a sophisticated SDP-based initialization for determining an initial set of centroids.

Table \ref{tab:heuristic} reports the frequency of each algorithm succeeding in finding the best solution among the competitors and the frequency of each algorithm succeeding in finding the solution obtained by PC-SOS-SDP on the 30 semi-supervised clustering instances of each dataset.
Note that, when the constraints set is feasible, both IPC-$k$-means and the post-processing procedure are guaranteed to return a valid clustering solution, while COP-$k$-means can potentially fail in finding any solution at all. The rate at which this event occurred in our experiments is reported in the rightmost COP-$k$-means column.
The results presented in Table \ref{tab:heuristic} show that IPC-$k$-means is almost unbeaten from any of the other heuristic algorithms, and it also succeeds in finding the optimal clustering solution on $79\%$ of our instances.

\begin{table}
\centering
\begin{tabular}{c|cc|cc|cc|ccc}
\toprule
             &    \multicolumn{2}{c|}{IPC-$k$-means} & \multicolumn{2}{c|}{BLPKM$^{++}$} & \multicolumn{2}{c|}{PP} & \multicolumn{3}{c}{COP-$k$-means} \\
            \cmidrule(lr){2-10}
            dataset & \multicolumn{1}{c|}{best} & opt & \multicolumn{1}{c|}{best} & opt & \multicolumn{1}{c|}{best} & opt & \multicolumn{1}{c|}{best} & \multicolumn{1}{c|}{opt} & fail\\
\midrule
               Iris &  96.7\% &  96.7\% &           60\% &    60\% &       16.7\% & 16.7\% &       3.3\% &  3.3\% & 6.7\%    \\
               Wine &   100\% &   100\% &         56.7\% &  56.7\% &       16.7\% & 16.7\% &       3.3\% &  3.3\% & 3.3\%    \\
              Seeds &  96.7\% &  96.7\% &           70\% &  66.7\% &         30\% &   30\% &      16.7\% & 16.7\% & 20\%     \\
      Connectionist &  96.7\% &  93.3\% &         26.7\% &  26.7\% &        6.7\% &  6.7\% &         0\% &    0\% & 63.3\%   \\
              Glass &   100\% &  43.3\% &            0\% &     0\% &          0\% &    0\% &         0\% &    0\% &          \\
              Heart &   100\% &  96.7\% &           10\% &    10\% &          0\% &    0\% &         0\% &    0\% & 63.3\%   \\
          Vertebral &  96.7\% &  96.7\% &         73.3\% &  73.3\% &          0\% &    0\% &         0\% &    0\% & 66.7\%   \\
             Accent &  96.7\% &  33.3\% &            0\% &     0\% &        3.3\% &    0\% &         0\% &    0\% &          \\
              Ecoli &   100\% &    40\% &            0\% &     0\% &          0\% &    0\% &         0\% &    0\% &          \\
            ECG5000 &   100\% &  63.3\% &            0\% &     0\% &          0\% &    0\% &         0\% &    0\% &          \\
          Computers &   100\% &  96.7\% &         16.7\% &  13.3\% &          0\% &    0\% &         0\% &    0\% & 26.7\%   \\
               Gene &  96.7\% &  96.7\% &           50\% &  46.7\% &       43.3\% &   40\% &      46.7\% & 46.7\% &          \\
\bottomrule
\end{tabular}
\caption{Rates of being the best among the considered heuristic algorithms (on the left) and of finding the optimal clustering (on the right); for COP-$k$-means the frequency of not being able to produce any solution at all is reported in the rightmost column}
\label{tab:heuristic}
\end{table}

\begin{figure}[!ht]
\label{heuristic_boxplot}
\centering
\includegraphics[scale=1]{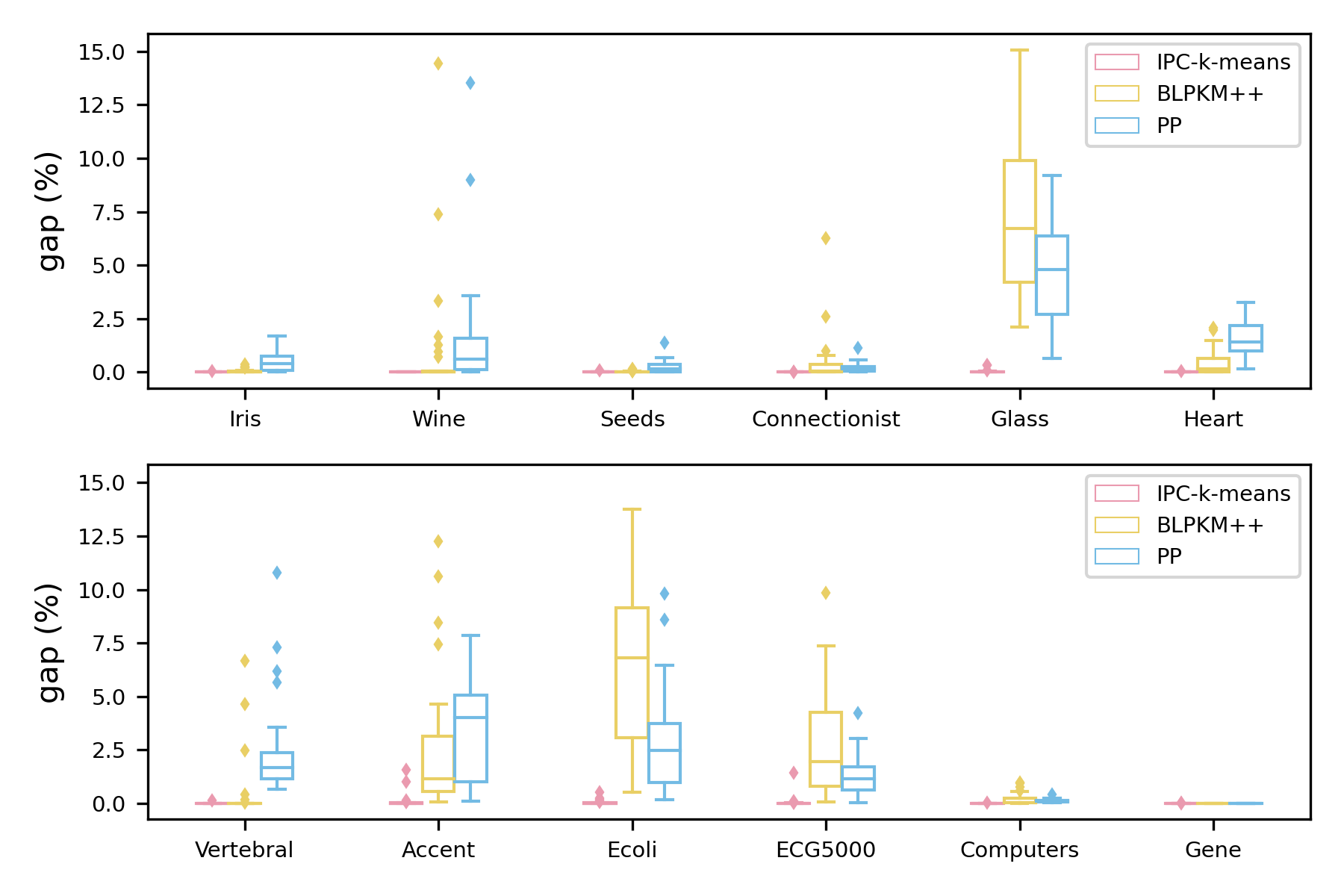}
\caption{Relative gaps between the heuristic solutions and the optimal clustering.}
\label{fig:heuristic_gap}
\end{figure}

The boxplot in Figure \ref{fig:heuristic_gap} shows a comparison among the relative gaps obtained by each of the four algorithms and the optimal solution on all the 30 instances of each dataset.

The exact computational time required for running IPC-$k$-means on each instance, including the time spent for the SDP-based initialization, can be found in \ref{app:num} since it coincides with the running time required for processing the root node in PC-SOS-SDP.
Naturally, relying on the resolution of an SDP problem, our heuristic does not scale well with the dataset cardinality. However, we have shown that it is able to produce better solutions in comparison to the other heuristic algorithms that we considered. Considering the scalability of SDPNAL+ and Gurobi, this heuristic can still be efficiently used on instances consisting of up to a few thousands data points. Indeed, solving the SDP (\ref{prob:SDPbab}) to moderate precision with data points ranging from $1000$ to $2000$ and a few thousands of inequalities {takes between 10 and 30 minutes} \citep{yang2015sdpnal}, whereas the time needed by Gurobi for the assignment problem (\ref{prob:ass}) is much smaller \citep{gurobi}. Clearly, the exact algorithm can not be applied on larger instances since a very large number of SDPs may be needed to be solved, making the approach impractical.

\subsection{Clustering evaluation}
Evaluating the quality of the solution produced by clustering algorithms is a difficult problem in general, since there is no ``gold standard'' solution in cluster analysis.  In addition to the MSSC objective function, we consider two widely used metrics for cluster evaluation.
Following the semi-supervised literature on this subject, we measure the agreement between the produced clustering solution and the true solution which is obtained on the basis of given class labels. To this end, the Adjusted Rand Index (ARI) \citep{hubert1985comparing} and the Adjusted Mutual Information (AMI) \citep{vinh2010information} are recommended as the indices of choice for measuring agreement between two partitions in semi-supervised clustering analysis. 
Denote by $U = \{C_1, \dots, C_k\}$ the ground truth class assignment and by $V = \{\hat{C}_1, \dots, \hat{C}_k\}$ a clustering of $n$ data points. The overlap between partitions $U$ and $V$ can be summarized in a contingency table where each entry $c_{ij}$ represents the number of data points that are common to clusters $C_i$ and $\hat{C}_j$, i.e., $c_{ij} = |C_i \cap \hat{C}_j |$. Furthermore, let $a_i$ and $b_j$ be the partial sums of the contingency table, i.e., $a_i = \sum_{j=1}^{k} c_{ij}$ and $b_j = \sum_{i=1}^{k} c_{ij}$. The ARI is computed as
\begin{equation*}
    ARI(U, V) = \frac{\sum_{i=1}^k \sum_{j=1}^k \binom{c_{ij}}{2} - \big[\sum_{i=1}^{k} \binom{a_i}{2} \sum_{j=1}^{k} \binom{b_j}{2}\big] / \binom{n}{2}}{\frac{1}{2}\big[\sum_{i=1}^{k} \binom{a_i}{2} \sum_{j=1}^{k} \binom{b_j}{2}\big] - \big[\sum_{i=1}^{k} \binom{a_i}{2} \sum_{j=1}^{k} \binom{b_j}{2}\big] / \binom{n}{2}}.
\end{equation*}
It takes a value of 1 when the two partitions are identical, the value 0 when the index equals its expected value and it can yield negative values if the index is less than the expected index.
The AMI is an information theoretic measure that quantifies the information shared by the partitions. It is computed as
\begin{equation*}
    AMI(U, V) = \frac{MI(U, V) - \mathbb{E}[MI(U, V)]}{\max(H(U), H(V)) - \mathbb{E}[MI(U, V)]},
\end{equation*}
where $MI(U,V)$ and $\mathbb{E}[MI(U, V)]$ are the mutual information and the expected mutual information between $U$ and $V$, and $H(U)$, $H(V)$ are the entropy associated with $U$ and $V$, respectively. These quantities are defined as
\begin{equation*}
    MI(U, V) = \sum_{i=1}^{k} \sum_{j=1}^{k} \frac{c_{ij}}{n} \textrm{log}\bigg(\frac{n c_{ij}}{a_i b_j}\bigg), \quad H(U) = -\sum_{i=1}^{k} \frac{a_i}{n} \textrm{log}\bigg(\frac{a_i}{n}\bigg), \quad H(V) = -\sum_{j=1}^{k} \frac{b_j}{n} \textrm{log}\bigg(\frac{b_j}{n}\bigg),
\end{equation*}
\begin{equation*}
    \mathbb{E}[MI(U, V)] = \sum_{i=1}^{k} \sum_{j=1}^{k} \sum_{n_{ij} = \max(1, a_i+b_j-n)}^{\min(a_i, b_j)} \frac{c_{ij}}{n} \textrm{log}\bigg( \frac{n c_{ij}}{a_i b_j}\bigg) \frac{a_i!b_j!(n-a_i)!(n-b_j)!}{n!c_{ij}!(a_i - c_{ij})!(b_j - c_{ij})!(n-a_i-b_j+c_{ij})!}.
\end{equation*}
The AMI takes a value of 1 when the two partitions perfectly match and 0 when the mutual information
between the two partitions equals its expected value. Note that the AMI can take negative values since random partitions (independent labellings) have an expected AMI around 0 on average.
\\
In this section, we address the following research questions: 
\begin{enumerate}
 \item Does the cluster quality increase when ML and/or CL constraints are added?
 \item Is there an advantage in applying an exact algorithm, i.e., does a more accurate solution of the semi-supervised MSSC correspond to a higher cluster quality?
\end{enumerate}

In Tables \ref{fig:metrics_ml}, \ref{fig:metrics_cl} and \ref{fig:metrics_mlcl} for each dataset we consider the instances where only ML, only CL and both ML and CL constraints are included, respectively. In all these tables, we report:
 \begin{itemize}
  \item violated: the average percentage of pairwise constraints that are violated in the global minimum of the unconstrained MSSC.
  \item MSSC ARI (AMI): the ARI (AMI) computed on the global solution of the unconstrained MSSC (by means of the exact algorithm proposed in \cite{piccialli2021sos}).
  \item PP ARI (AMI): the ARI (AMI) computed on the feasible clustering produced by the heuristic algorithm in \cite{vrain2020constrained}.
  \item PC-SOS-SDP ARI (AMI): the ARI (AMI) computed on the solution produced by the exact algorithm proposed in this paper.
 \end{itemize}
 Looking at Tables \ref{fig:metrics_ml}-\ref{fig:metrics_mlcl}, it emerges that there is always an increase in both indices (ARI and AMI) when any pairwise constraint is added. This observation holds even if the solution is not optimal (the average ARI and AMI of the PP solution are already higher than the ones of the unconstrained MSSC). However,  when the exact solver finds a better solution, there is a further improvement in the indices. In some cases (see, for example, the Heart dataset), the starting value of the ARI/AMI indices is close to zero, meaning that the unconstrained solution does not overlap with the ground truth. For these examples, the indices increase when ML and CL are included, but they remain low. When this is the case, a kernel-based clustering (see \cite{filippone2008survey}) may be a viable approach. In some cases, there is already a very good agreement between the unconstrained solution and the semi-supervised one (see, for example, the Gene dataset). Also, in this case, the indices increase and get very close to 1 when constraints are added. Note that there is some correlation between the percentage of violated constraints and the quality of the unconstrained solution. Indeed, if the agreement between the ground truth and the unconstrained solution is low, most of the randomly generated constraints will be violated. When there is a good agreement, it becomes harder to randomly generate constraints violated by the unconstrained solution. Therefore, the constraints become less informative. There are some intermediate cases, where adding the constraints leads to a significant improvement in the solution quality: see, for example, the datasets Wine, Ecoli, and ECG5000. Not surprisingly, the most informative constraints are the must link, allowing a higher improvement in both ARI and AMI. Overall, these results confirm the importance of using background knowledge, if available, since they show the increased quality of the obtained clustering. Furthermore, finding the globally optimal constrained clustering solution always translates into better clustering results.

\begin{table}[!ht]
    \centering
    \begin{tabular}{cc|ccc|ccc}
\toprule
      & &    \multicolumn{3}{c|}{ARI} & \multicolumn{3}{c}{AMI} \\
            \cmidrule(lr){3-8} 
            dataset & violated & \multicolumn{1}{c|}{MSSC} &  \multicolumn{1}{c|}{PP} & PC-SOS-SDP & \multicolumn{1}{c|}{MSSC} & \multicolumn{1}{c|}{PP} & PC-SOS-SDP\\
\midrule
         Iris & 16.1\% &     0.71 & 0.874 & 0.891 &    0.741 & 0.871 & 0.886 \\
         Wine & 41.3\% &    0.361 & 0.519 & 0.528 &    0.413 & 0.538 & 0.563 \\
Connectionist & 51.6\% &    0.005 & 0.051 & 0.065 &    0.003 & 0.037 & 0.056 \\
        Seeds & 18.4\% &    0.707 & 0.822 & 0.847 &    0.682 & 0.790 & 0.810 \\
        Glass & 57.7\% &    0.161 & 0.211 & 0.269 &    0.294 & 0.329 & 0.394 \\
        Heart & 46.0\% &   -0.004 & 0.015 & 0.464 &   -0.003 & 0.003 & 0.350 \\
    Vertebral & 42.1\% &    0.105 & 0.067 & 0.381 &    0.253 & 0.233 & 0.389 \\
       Accent & 76.7\% &    0.032 & 0.211 & 0.291 &    0.162 & 0.237 & 0.350 \\
        Ecoli & 57.9\% &    0.426 & 0.569 & 0.731 &    0.555 & 0.671 & 0.736 \\
      ECG5000 & 50.0\% &    0.465 & 0.582 & 0.611 &     0.51 & 0.576 & 0.603 \\
    Computers & 46.6\% &   -0.001 & 0.025 & 0.038 &   -0.001 & 0.019 & 0.035 \\
         Gene & 1.3\%  &    0.971 & 0.993 & 0.998 &    0.962 & 0.987 & 0.993 \\
\bottomrule
\end{tabular}
\caption{Clustering evaluation metrics computed on solutions of instances with only ML constraints.}
\label{fig:metrics_ml}
\end{table}

\begin{table}[!ht]
\centering
\begin{tabular}{cc|ccc|ccc}
\toprule
      & &    \multicolumn{3}{c|}{ARI} & \multicolumn{3}{c}{AMI} \\
            \cmidrule(lr){3-8} 
            dataset & violated & \multicolumn{1}{c|}{MSSC} &  \multicolumn{1}{c|}{PP} & PC-SOS-SDP & \multicolumn{1}{c|}{MSSC} & \multicolumn{1}{c|}{PP} & PC-SOS-SDP\\
\midrule
         Iris &  9.7\%  &    0.71 & 0.816 & 0.836 &    0.741 & 0.815 & 0.829 \\
         Wine &  20.8\% &   0.361 & 0.409 & 0.429 &    0.413 & 0.436 & 0.454 \\
Connectionist &  51.6\% &   0.005 & 0.056 & 0.077 &    0.003 & 0.041 & 0.061 \\
        Seeds &  10.7\% &   0.707 & 0.792 & 0.812 &    0.682 & 0.764 & 0.779 \\
        Glass &  24.0\% &   0.161 & 0.169 & 0.205 &    0.294 & 0.290 & 0.322 \\
        Heart &  53.4\% &  -0.004 & 0.022 & 0.463 &   -0.003 & 0.015 & 0.376 \\
    Vertebral &  47.0\% &   0.105 & 0.273 & 0.375 &    0.253 & 0.279 & 0.345 \\
       Accent &  18.7\% &   0.032 & 0.028 & 0.057 &    0.162 & 0.180 & 0.203 \\
        Ecoli &  3.4\%  &   0.426 & 0.449 & 0.452 &    0.555 & 0.592 & 0.613 \\
      ECG5000 &  5.2\%  &   0.465 & 0.459 & 0.512 &     0.51 & 0.522 & 0.553 \\
    Computers &  52.3\% &  -0.001 & 0.013 & 0.041 &   -0.001 & 0.010 & 0.037 \\
         Gene &  0.5\%  &   0.971 & 0.989 & 0.995 &    0.962 & 0.981 & 0.986 \\
\bottomrule
\end{tabular}
\caption{Clustering evaluation metrics computed on solutions of instances with only CL constraints.}
\label{fig:metrics_cl}
\end{table}

\begin{table}[!ht]
\centering
\begin{tabular}{cc|ccc|ccc}
\toprule
      & &    \multicolumn{3}{c|}{ARI} & \multicolumn{3}{c}{AMI} \\
            \cmidrule(lr){3-8} 
            dataset & violated & \multicolumn{1}{c|}{MSSC} &  \multicolumn{1}{c|}{PP} & PC-SOS-SDP & \multicolumn{1}{c|}{MSSC} & \multicolumn{1}{c|}{PP} & PC-SOS-SDP\\
\midrule
         Iris & 14.9\% &     0.71 & 0.848 & 0.873 &    0.741 & 0.844 & 0.863 \\
         Wine & 32.8\% &    0.361 & 0.457 & 0.491 &    0.413 & 0.484 & 0.512 \\
Connectionist & 47.2\% &    0.005 & 0.042 & 0.088 &    0.003 & 0.031 & 0.070 \\
        Seeds & 16.0\% &    0.707 & 0.820 & 0.832 &    0.682 & 0.789 & 0.798 \\
        Glass & 40.9\% &    0.161 & 0.195 & 0.223 &    0.294 & 0.320 & 0.345 \\
        Heart & 50.4\% &   -0.004 & 0.005 & 0.447 &   -0.003 & 0.005 & 0.352 \\
    Vertebral & 44.2\% &    0.105 & 0.229 & 0.355 &    0.253 & 0.262 & 0.332 \\
       Accent & 46.3\% &    0.032 & 0.064 & 0.071 &    0.162 & 0.192 & 0.234 \\
        Ecoli & 31.9\% &    0.426 & 0.463 & 0.554 &    0.555 & 0.619 & 0.667 \\
      ECG5000 & 27.7\% &    0.465 & 0.484 & 0.489 &     0.51 & 0.534 & 0.555 \\
    Computers & 49.8\% &   -0.001 & 0.015 & 0.035 &   -0.001 & 0.011 & 0.032 \\
         Gene & 0.9\%  &    0.971 & 0.991 & 0.997 &    0.962 & 0.985 & 0.991 \\
\bottomrule
\end{tabular}
\caption{Clustering evaluation metrics computed on solutions of instances with both ML and CL constraints.}
\label{fig:metrics_mlcl}
\end{table}

\section{Conclusions}
In this paper, we have proposed PC-SOS-SDP: a branch-and-{cut} algorithm for semi-supervised MSSC with pairwise constraints. For the lower bound procedure, we use an SDP relaxation that exploits ML constraints to reduce the size of the problem. In addition, we add three types of valid inequalities in a cutting plane fashion to generate tight bounds. We have also developed a heuristic named IPC-$k$-means that is reminiscent of the popular $k$-means algorithm. When the problem is feasible this heuristic returns a good quality clustering and an upper bound on the optimality gap of the provided solution; otherwise it returns a certificate of infeasibility.
Numerical results of the overall branch-and-{cut} algorithm impressively exhibit the efficiency of PC-SOS-SDP: we can solve real-world instances up to 800 data points with different combinations of ML and CL constraints and with a huge number of features. To the best of our knowledge, no other exact methods can handle generic instances of that size. PC-SOS-SDP can deal with other constrained versions of MSSC like those with diameter and split constraints \citep{davidson2005clustering}. These constraints can be represented as a disjunction and conjunction of ML and CL constraints, thus making their implementation easy in our specialized solver. As future work, we plan to extend PC-SOS-SDP for semi-supervised MSSC with cluster-level constraints concerning the cardinality of the clusters. 

\newpage
\appendix
\section{Proof of Theorem 3.2}\label{app:proof}
\begin{proof}
We introduce dual variables $y \in \mathbb{R}^{s+1}$, $\mu \in \mathbb{R}^{|\mathcal{\overline{CL}}|}$, $\tilde{V} \in \mathcal{S}^s$ and we write the dual of Problem \eqref{prob:SDPbab} (omitting the constant part of the objective function) as
\begin{maxi!}[2]
{}{\sum_{i=1}^{s} y_i + k y_{s+1}}
    {\label{prob:dualSDPbab}}{}
\addConstraint{-\mathcal{T}^s W (\mathcal{T}^s)^\top - \sum_{i=1}^{s} y_i E_i^s - y_{s+1} \Diag(e^s) - \sum_{(i, j) \in \mathcal{\overline{CL}}} \mu_{\tau(i,j)} E_{ij} - S}{= V}{}
\addConstraint{S \in \mathcal{S}_+^s, \ V \geq 0}{}{}.
\end{maxi!}
Given $Z^\star$, first observe that the dual objective function can be written as 
\begin{equation*}
    \sum_{i=1}^{s} y_i + k y_{s+1} + \sum_{(i, j) \in \mathcal{\overline{CL}}} \mu_{\tau(i,j)} \inprod{E_{ij}}{Z^\star}
\end{equation*}
since for any feasibile $Z$ of the primal problem, $Z_{ij} = 0$ for $(i, j) \in \mathcal{\overline{CL}}$.
In order to show that $p^\star \ge lb$, consider the difference between the primal and the dual objective as
\begin{subequations}
\begin{align*}
    \inprod{-\mathcal{T}^s W (\mathcal{T}^s)^\top}{Z^\star} & - \Big(\sum_{i=1}^{s} \tilde{y}_i + k \tilde{y}_{s+1} + \sum_{(i, j) \in \mathcal{\overline{CL}}} \tilde{\mu}_{\tau(i,j)} \inprod{E_{ij}}{Z^\star} \Big) \\
    & = \inprod{-\mathcal{T}^s W (\mathcal{T}^s)^\top}{Z^\star} - \sum_{i=1}^{s} \tilde{y}_i \inprod{E_i^s}{Z^\star} - \inprod{\Diag(e^s)}{Z^\star} \tilde{y}_{s+1} - \sum_{(i, j) \in \mathcal{\overline{CL}}} \tilde{\mu}_{\tau(i,j)} \inprod{E_{ij}}{Z^\star} \\
    & = \inprod{-\mathcal{T}^s W (\mathcal{T}^s)^\top - \sum_{i=1}^{s} \tilde{y}_i E_i^s - \Diag(e^s) \tilde{y}_{s+1} - \sum_{(i, j) \in \mathcal{\overline{CL}}} \tilde{\mu}_{\tau(i,j)} E_{ij}}{Z^\star}\\
    & = \inprod{\tilde{V} + \tilde{S}}{Z^\star}.
\end{align*}
\end{subequations}
Using Lemma \ref{lemma:jansson}, we obtain
\begin{align*}
\inprod{\tilde{V}}{Z^\star} + \inprod{\tilde{S}}{Z^\star} \geq \inprod{\tilde{V}}{Z^\star} + \bar{z} \sum_{i\colon \lambda_i(\tilde{S}) < 0} \lambda_i(\tilde{S}) \geq \bar{z} \sum_{i\colon \lambda_i(\tilde{S}) < 0} \lambda_i(\tilde{S}),
\end{align*}
where the last inequality holds because $\tilde{V}$ is nonnegative.
\end{proof}
\newpage
\section{PC-SOS-SDP Detailed Numerical Results}\label{app:num}
Here we report a more detailed description of the resolution of each of the 360 instances generated for our computational experiments.
The instances are divided by datasets and presented in Tables \ref{tab:iris_detailed} - \ref{tab:gene_detailed}.
For each instance we report the following data:
\begin{itemize}
    \item \textit{ml}: number of \textit{must-link} constraints;
    \item \textit{cl}: number of \textit{cannot-link} constraints;
    \item \textit{seed}: seed used in the random generation process;
    \item \textit{size}: number of separate data points on the root node (it can differ from the dataset size in the presence of \textit{must-link} constraints);
    \item \textit{f}: the optimal solution value;
    \item \textit{cp\textsubscript{0}}: the number of cutting plane iterations performed on the root node for the lower bound computation;
    \item \textit{ineq\textsubscript{0}}: the number of inequalities of the last SDP solved at the root in the cutting-plane procedure;
    \item \textit{gap\textsubscript{0}}: the relative gap between the root lower and upper bound calculated as $(UB - LB) / UB$;
    \item \textit{time\textsubscript{0}}: running time on the root node, expressed in seconds;
    \item \textit{nodes}: number of branch-and-bound nodes explored;
    \item \textit{time}: total running time, expressed in seconds;

\end{itemize}

\begin{table}[!h]
\centering
\footnotesize
    \begin{tabular}{ccccccccccc}
\toprule
 ml &  cl &  seed &  size &       f &  cp\textsubscript{0} &  ineq\textsubscript{0} &   gap\textsubscript{0} & time\textsubscript{0}(s) & nodes &  time(s) \\
\midrule
 50 &   0 &     0 &   101 & 83.6299 &     2 &    6198 & 0.0004\% &          9 &      1 &       12 \\
 50 &   0 &     1 &   100 & 85.1940 &     1 &    5356 & 0.0010\% &          5 &      1 &        7 \\
 50 &   0 &     2 &   100 & 87.6941 &     1 &    5419 & 0.0007\% &          8 &      1 &       11 \\
 50 &   0 &     3 &   100 & 84.7522 &     1 &    5365 & 0.0060\% &          5 &      1 &        8 \\
 50 &   0 &     4 &   101 & 83.5486 &     2 &    5301 & 0.0100\% &         25 &      5 &       45 \\\hline
 25 &  25 &     0 &   125 & 79.9578 &     2 &    5742 & 0.0013\% &         21 &      1 &       23 \\
 25 &  25 &     1 &   125 & 86.7623 &     2 &    6149 & 0.0003\% &         18 &      1 &       21 \\
 25 &  25 &     2 &   125 & 82.0447 &     2 &    5854 & 0.0284\% &         12 &      3 &       35 \\
 25 &  25 &     3 &   125 & 81.4878 &     3 &    7231 & 0.0013\% &         35 &      1 &       38 \\
 25 &  25 &     4 &   125 & 81.5452 &     2 &    5717 & 0.0363\% &          8 &      3 &       24 \\\hline
  0 &  50 &     0 &   150 & 84.1465 &     3 &    9887 & 0.0047\% &         77 &      1 &       80 \\
  0 &  50 &     1 &   150 & 85.6707 &     5 &    3779 & 0.0172\% &         93 &      5 &      178 \\
  0 &  50 &     2 &   150 & 80.2087 &     4 &    6353 & 0.0173\% &         53 &      3 &       74 \\
  0 &  50 &     3 &   150 & 84.8764 &     2 &    6362 & 0.0012\% &         29 &      1 &       31 \\
  0 &  50 &     4 &   150 & 85.0928 &     2 &    5686 & 0.0103\% &         26 &      1 &       29 \\\hline
100 &   0 &     0 &    57 & 85.6052 &     1 &    3732 & 0.0011\% &          3 &      1 &        6 \\
100 &   0 &     1 &    55 & 87.9862 &     1 &    5093 & 0.0007\% &          4 &      1 &        6 \\
100 &   0 &     2 &    54 & 87.9577 &     2 &    2709 & 0.0026\% &          3 &      1 &        5 \\
100 &   0 &     3 &    59 & 84.8172 &     1 &    5136 & 0.0004\% &          3 &      1 &        6 \\
100 &   0 &     4 &    55 & 87.0724 &     1 &    3395 & 0.0001\% &          3 &      1 &        6 \\\hline
 50 &  50 &     0 &   100 & 84.5632 &     2 &    5311 & 0.0132\% &         11 &      3 &       23 \\
 50 &  50 &     1 &   100 & 83.3805 &     1 &    5352 & 0.0012\% &          3 &      1 &        6 \\
 50 &  50 &     2 &   100 & 86.8824 &     2 &    5526 & 0.0052\% &         13 &      1 &       15 \\
 50 &  50 &     3 &   100 & 84.9081 &     1 &    5299 & 0.0003\% &          3 &      1 &        5 \\
 50 &  50 &     4 &   100 & 88.1519 &     1 &    5237 & 0.0008\% &          4 &      1 &        7 \\\hline
  0 & 100 &     0 &   150 & 81.5149 &     3 &    6156 & 0.0031\% &         43 &      1 &       45 \\
  0 & 100 &     1 &   150 & 83.2654 &     2 &    6017 & 0.0716\% &         26 &      3 &       46 \\
  0 & 100 &     2 &   150 & 86.1146 &     3 &    8540 & 0.0248\% &         50 &      3 &       76 \\
  0 & 100 &     3 &   150 & 87.2248 &     2 &    6391 & 0.3634\% &         40 &      3 &       75 \\
  0 & 100 &     4 &   150 & 86.6330 &     2 &    7546 & 0.0495\% &         25 &      3 &       51 \\
\bottomrule
\end{tabular}

    \caption{Dataset \textit{Iris}, $k$=3.}
    \label{tab:iris_detailed}
\end{table}

\begin{table}[!h]
\centering
\footnotesize
    \begin{tabular}{ccccccccccc}
\toprule
 ml &  cl &  seed &  size &         f &  cp\textsubscript{0} &  ineq\textsubscript{0} &   gap\textsubscript{0} & time\textsubscript{0}(s) &  nodes &  time(s) \\
\midrule
 50 &   0 &     0 &   128 & 3.62636e+06 &     4 &    4286 & 0.0117\% &         49 &      1 &       51 \\
 50 &   0 &     1 &   128 & 3.64084e+06 &     6 &    8179 & 0.0074\% &        138 &      1 &      141 \\
 50 &   0 &     2 &   128 & 3.33862e+06 &     3 &    9597 & 0.0041\% &         48 &      1 &       50 \\
 50 &   0 &     3 &   128 & 3.30173e+06 &     6 &    8268 & 0.0331\% &         95 &      3 &      118 \\
 50 &   0 &     4 &   129 & 3.51784e+06 &     5 &   10713 & 0.0233\% &        117 &      3 &      149 \\\hline
 25 &  25 &     0 &   153 & 3.23052e+06 &     3 &    6765 & 0.0299\% &         40 &      3 &       61 \\
 25 &  25 &     1 &   153 & 3.13854e+06 &     7 &    6174 & 0.0001\% &         98 &      1 &      100 \\
 25 &  25 &     2 &   153 & 3.28121e+06 &     3 &    3223 & 0.0041\% &         56 &      1 &       58 \\
 25 &  25 &     3 &   153 & 3.14043e+06 &     3 &    6403 & 0.0235\% &         48 &      1 &       50 \\
 25 &  25 &     4 &   153 & 2.95302e+06 &     4 &    4415 & 0.0209\% &         57 &      1 &       60 \\\hline
  0 &  50 &     0 &   178 & 2.60234e+06 &     3 &    7192 & 0.0142\% &         42 &      1 &       44 \\
  0 &  50 &     1 &   178 & 2.86950e+06 &     4 &    6385 & 0.0014\% &         81 &      1 &       83 \\
  0 &  50 &     2 &   178 & 2.63343e+06 &     3 &    7186 & 0.1571\% &         54 &      3 &       78 \\
  0 &  50 &     3 &   178 & 2.97316e+06 &     4 &    6263 & 0.0028\% &         87 &      1 &       89 \\
  0 &  50 &     4 &   178 & 2.73127e+06 &     5 &    6259 & 0.0174\% &        129 &      1 &      132 \\\hline
100 &   0 &     0 &    80 & 4.19548e+06 &     1 &    5671 & 0.0099\% &          4 &      1 &        6 \\
100 &   0 &     1 &    80 & 4.43060e+06 &     2 &   10046 & 0.0018\% &          8 &      1 &       10 \\
100 &   0 &     2 &    83 & 4.30206e+06 &     4 &    8660 & 0.0334\% &         86 &      3 &      108 \\
100 &   0 &     3 &    82 & 4.04948e+06 &     3 &    7119 & 0.0007\% &         22 &      1 &       24 \\
100 &   0 &     4 &    82 & 4.02972e+06 &     1 &    5563 & 0.0004\% &          7 &      1 &       10 \\\hline
 50 &  50 &     0 &   129 & 3.56553e+06 &     5 &    5844 & 0.0775\% &        128 &      3 &      151 \\
 50 &  50 &     1 &   129 & 3.75127e+06 &     5 &    8668 & 0.0264\% &         84 &      5 &      174 \\
 50 &  50 &     2 &   128 & 3.74904e+06 &     2 &    6734 & 0.0021\% &         21 &      1 &       23 \\
 50 &  50 &     3 &   128 & 4.03569e+06 &     3 &    8517 & 0.0025\% &         59 &      1 &       61 \\
 50 &  50 &     4 &   128 & 4.01734e+06 &     3 &    6883 & 0.0001\% &         53 &      1 &       56 \\\hline
  0 & 100 &     0 &   178 & 3.09635e+06 &     4 &    8677 & 0.0125\% &         86 &      1 &       89 \\
  0 & 100 &     1 &   178 & 3.09061e+06 &     6 &    6676 & 0.0190\% &        156 &      3 &      196 \\
  0 & 100 &     2 &   178 & 3.12229e+06 &     5 &    6284 & 0.0393\% &        142 &      9 &      270 \\
  0 & 100 &     3 &   178 & 3.33813e+06 &     5 &    6330 & 0.2200\% &        101 &     11 &      275 \\
  0 & 100 &     4 &   178 & 2.90243e+06 &     4 &    6420 & 0.0221\% &         92 &      7 &      253 \\
\bottomrule
\end{tabular}

    \caption{Dataset \textit{Wine}, $k$=3.}
    \label{tab:wine_detailed}
\end{table}

\begin{table}[!h]
\centering
\footnotesize
    \begin{tabular}{ccccccccccc}
\toprule
 ml &  cl &  seed &  size &       f &  cp\textsubscript{0} &  ineq\textsubscript{0} &   gap\textsubscript{0} & time\textsubscript{0}(s) &  nodes &  time(s) \\
\midrule
 50 &   0 &     0 &   158 & 296.063 &     4 &    5703 & 0.0050\% &        294 &      1 &      297 \\
 50 &   0 &     1 &   159 & 301.857 &     6 &    5995 & 0.0010\% &        113 &      1 &      116 \\
 50 &   0 &     2 &   158 & 299.851 &     5 &    6085 & 0.0062\% &        180 &      1 &      182 \\
 50 &   0 &     3 &   158 & 304.610 &     5 &    5918 & 0.0040\% &        100 &      1 &      102 \\
 50 &   0 &     4 &   158 & 304.781 &     6 &    6348 & 0.0066\% &        104 &      1 &      107 \\\hline
 25 &  25 &     0 &   176 & 298.022 &     6 &    5142 & 0.0053\% &        113 &      1 &      116 \\
 25 &  25 &     1 &   181 & 302.256 &     5 &    6001 & 0.0060\% &        133 &      1 &      135 \\
 25 &  25 &     2 &   177 & 305.470 &    10 &    5965 & 0.0023\% &        352 &      1 &      355 \\
 25 &  25 &     3 &   177 & 305.012 &     5 &    6078 & 0.0039\% &        138 &      1 &      140 \\
 25 &  25 &     4 &   177 & 300.234 &     7 &    6190 & 0.0055\% &        377 &      1 &      380 \\\hline
  0 &  50 &     0 &   189 & 302.604 &     5 &    6147 & 0.0003\% &        241 &      1 &      243 \\
  0 &  50 &     1 &   189 & 299.810 &     5 &    6438 & 0.0006\% &        110 &      1 &      113 \\
  0 &  50 &     2 &   182 & 301.590 &     8 &    5965 & 0.0062\% &        192 &      1 &      195 \\
  0 &  50 &     3 &   190 & 300.461 &     5 &    6414 & 0.0002\% &        118 &      1 &      120 \\
  0 &  50 &     4 &   183 & 308.610 &     6 &    5243 & 0.0045\% &        154 &      1 &      157 \\\hline
100 &   0 &     0 &   108 & 331.031 &     5 &    5805 & 0.0131\% &         49 &      3 &       68 \\
100 &   0 &     1 &   108 & 327.204 &     4 &    5618 & 0.0059\% &         49 &      1 &       52 \\
100 &   0 &     2 &   109 & 322.181 &     6 &    6458 & 0.0015\% &         92 &      1 &       95 \\
100 &   0 &     3 &   110 & 326.708 &     6 &    5614 & 0.0161\% &         71 &      3 &      108 \\
100 &   0 &     4 &   110 & 331.773 &     5 &    5716 & 0.0009\% &         52 &      1 &       55 \\\hline
 50 &  50 &     0 &   135 & 329.379 &     5 &    1932 & 0.0013\% &         99 &      1 &      101 \\
 50 &  50 &     1 &   139 & 330.708 &     6 &    3704 & 0.0044\% &        117 &      1 &      119 \\
 50 &  50 &     2 &   144 & 324.190 &     4 &    5820 & 0.0002\% &         84 &      1 &       86 \\
 50 &  50 &     3 &   141 & 329.774 &     3 &    5637 & 0.0041\% &         42 &      1 &       45 \\
 50 &  50 &     4 &   139 & 322.186 &     6 &    3629 & 0.0055\% &         86 &      1 &       89 \\\hline
  0 & 100 &     0 &   148 & 329.864 &     6 &    4583 & 0.0049\% &        137 &      1 &      139 \\
  0 & 100 &     1 &   139 & 325.496 &     4 &    7230 & 0.0018\% &        108 &      1 &      110 \\
  0 & 100 &     2 &   133 & 316.950 &     4 &    5555 & 0.0069\% &         61 &      1 &       64 \\
  0 & 100 &     3 &   137 & 331.657 &     3 &    5661 & 0.0001\% &         61 &      1 &       64 \\
  0 & 100 &     4 &   140 & 332.885 &     4 &    5672 & 0.0156\% &         74 &      3 &      136 \\
\bottomrule
\end{tabular}

    \caption{Dataset \textit{Connectionist}, $k$=2.}
    \label{tab:connectionist_detailed}
\end{table}

\begin{table}[!h]
\centering
\footnotesize
    \begin{tabular}{ccccccccccc}
\toprule
 ml &  cl &  seed &  size &       f &  cp\textsubscript{0} &  ineq\textsubscript{0} &   gap\textsubscript{0} & time\textsubscript{0}(s) &  nodes &  time(s) \\
\midrule
 50 &   0 &     0 &   160 & 620.233 &     3 &    5623 & 0.0005\% &         47 &      1 &       49 \\
 50 &   0 &     1 &   160 & 617.303 &     2 &    5834 & 0.0041\% &         29 &      1 &       32 \\
 50 &   0 &     2 &   160 & 619.150 &     4 &    5861 & 0.0193\% &         61 &      5 &      102 \\
 50 &   0 &     3 &   160 & 624.169 &     4 &    5925 & 0.0352\% &         64 &     13 &      194 \\
 50 &   0 &     4 &   160 & 616.667 &     2 &    5903 & 0.0989\% &         20 &      3 &       45 \\\hline
 25 &  25 &     0 &   185 & 616.182 &     9 &    5453 & 0.0408\% &        160 &     11 &      311 \\
 25 &  25 &     1 &   185 & 596.609 &     3 &    5923 & 0.0540\% &         46 &      5 &      119 \\
 25 &  25 &     2 &   185 & 617.730 &     4 &    9109 & 0.0492\% &         71 &      5 &      161 \\
 25 &  25 &     3 &   185 & 600.371 &     4 &    5969 & 0.0067\% &         63 &      1 &       66 \\
 25 &  25 &     4 &   185 & 608.407 &     5 &    5661 & 0.0194\% &         83 &      7 &      166 \\\hline
  0 &  50 &     0 &   210 & 601.960 &     4 &    6473 & 0.1434\% &         76 &     15 &      255 \\
  0 &  50 &     1 &   210 & 626.553 &     8 &    3577 & 0.0119\% &        152 &      3 &      186 \\
  0 &  50 &     2 &   210 & 605.245 &     6 &    6274 & 0.0498\% &        104 &     11 &      380 \\
  0 &  50 &     3 &   210 & 604.747 &     5 &    5972 & 0.0620\% &        103 &     15 &      213 \\
  0 &  50 &     4 &   210 & 603.729 &     3 &    7102 & 0.0032\% &        100 &      1 &      102 \\\hline
100 &   0 &     0 &   112 & 636.851 &     2 &    5796 & 0.0012\% &         17 &      1 &       19 \\
100 &   0 &     1 &   110 & 632.934 &     2 &    5537 & 0.0385\% &         15 &      5 &       39 \\
100 &   0 &     2 &   110 & 644.974 &     1 &    5299 & 0.0003\% &          8 &      1 &       11 \\
100 &   0 &     3 &   112 & 651.855 &     1 &    5422 & 0.0045\% &         12 &      1 &       15 \\
100 &   0 &     4 &   114 & 652.323 &     2 &    4820 & 0.0004\% &         10 &      1 &       12 \\\hline
 50 &  50 &     0 &   160 & 634.578 &     5 &    7792 & 0.0538\% &         68 &      7 &      146 \\
 50 &  50 &     1 &   160 & 641.830 &     2 &    6036 & 0.1153\% &         33 &      7 &      108 \\
 50 &  50 &     2 &   160 & 634.669 &     2 &    5690 & 0.0317\% &         28 &      3 &       57 \\
 50 &  50 &     3 &   160 & 641.844 &     2 &   10555 & 0.0004\% &         19 &      1 &       22 \\
 50 &  50 &     4 &   160 & 635.987 &     2 &    5840 & 0.0069\% &         25 &      1 &       28 \\\hline
  0 & 100 &     0 &   210 & 616.050 &     4 &    6241 & 0.1831\% &         73 &     33 &      206 \\
  0 & 100 &     1 &   210 & 611.600 &     4 &    5732 & 0.0731\% &         73 &      7 &      215 \\
  0 & 100 &     2 &   210 & 624.302 &     4 &    6436 & 0.1108\% &         89 &      3 &      131 \\
  0 & 100 &     3 &   210 & 656.197 &     7 &    6763 & 0.0917\% &        138 &      9 &      391 \\
  0 & 100 &     4 &   210 & 632.656 &     6 &    5856 & 0.0263\% &        130 &      3 &      186 \\
\bottomrule
\end{tabular}

    \caption{Dataset \textit{Seeds}, $k$=3.}
    \label{tab:seeds_detailed}
\end{table}

\begin{table}[!h]
\centering
\footnotesize
    \begin{tabular}{ccccccccccc}
\toprule
 ml &  cl &  seed &  size &       f &  cp\textsubscript{0} &  ineq\textsubscript{0} &   gap\textsubscript{0} & time\textsubscript{0}(s) &  nodes &  time(s) \\
\midrule
 50 &   0 &     0 &   165 &  86.1435 &     4 &    5787 & 0.0097\% &        104 &      1 &      107 \\
 50 &   0 &     1 &   164 &  94.1670 &     7 &    7296 & 0.0780\% &        121 &     23 &      377 \\
 50 &   0 &     2 &   164 &  98.5834 &     8 &    3789 & 0.0212\% &        125 &      3 &      148 \\
 50 &   0 &     3 &   165 & 103.6630 &     7 &    6280 & 0.0171\% &        176 &      3 &      208 \\
 50 &   0 &     4 &   164 & 100.8140 &     8 &    6767 & 0.0633\% &        196 &      5 &      331 \\\hline
 25 &  25 &     0 &   189 &  85.4317 &    11 &    2671 & 0.3518\% &        162 &     47 &      497 \\
 25 &  25 &     1 &   189 &  85.0901 &     8 &    6455 & 0.0827\% &        192 &     21 &      647 \\
 25 &  25 &     2 &   189 &  82.9471 &     6 &    6609 & 0.0004\% &        157 &      1 &      159 \\
 25 &  25 &     3 &   189 &  79.0361 &     9 &    2644 & 0.1926\% &        162 &    141 &      761 \\
 25 &  25 &     4 &   189 &  87.6976 &    13 &    2510 & 0.2106\% &        178 &     51 &      646 \\\hline
  0 &  50 &     0 &   214 &  76.2862 &    12 &    3211 & 0.1659\% &        222 &    115 &      720 \\
  0 &  50 &     1 &   214 &  75.4997 &    11 &    2942 & 0.0820\% &        237 &     13 &      416 \\
  0 &  50 &     2 &   214 &  76.0591 &    11 &    3595 & 0.0787\% &        218 &     11 &      455 \\
  0 &  50 &     3 &   214 &  75.8723 &    13 &    6813 & 0.1417\% &        301 &     15 &      731 \\
  0 &  50 &     4 &   214 &  76.8616 &    11 &    5588 & 0.0500\% &        204 &      9 &      375 \\\hline
100 &   0 &     0 &   120 & 105.7710 &     5 &    4089 & 0.0518\% &        100 &      5 &      219 \\
100 &   0 &     1 &   115 & 109.2680 &     4 &    6000 & 0.0074\% &        255 &      1 &      258 \\
100 &   0 &     2 &   119 & 113.0420 &     5 &    6618 & 0.0256\% &        307 &      3 &      403 \\
100 &   0 &     3 &   116 & 116.8250 &     5 &    6330 & 0.0805\% &        237 &      3 &      339 \\
100 &   0 &     4 &   119 & 107.3600 &     4 &    6296 & 0.1300\% &         70 &     11 &      343 \\\hline
 50 &  50 &     0 &   165 & 100.7680 &     9 &    3672 & 0.2211\% &        112 &     73 &      529 \\
 50 &  50 &     1 &   165 &  95.7760 &     8 &    4572 & 0.7826\% &        130 &     85 &      555 \\
 50 &  50 &     2 &   164 &  98.9496 &    10 &    4906 & 0.6066\% &        149 &     11 &      301 \\
 50 &  50 &     3 &   164 & 100.3210 &     6 &    6059 & 0.2626\% &        113 &     51 &      499 \\
 50 &  50 &     4 &   165 & 101.1380 &    10 &    2946 & 0.2828\% &         91 &     47 &      389 \\\hline
  0 & 100 &     0 &   214 &  76.7631 &     8 &    2620 & 0.1253\% &        178 &     11 &      328 \\
  0 & 100 &     1 &   214 &  77.3701 &    10 &    3643 & 0.2124\% &        230 &     81 &      773 \\
  0 & 100 &     2 &   214 &  78.1937 &     8 &    6346 & 0.0062\% &        159 &      1 &      162 \\
  0 & 100 &     3 &   214 &  76.7243 &    11 &    3811 & 0.1442\% &        198 &     21 &      527 \\
  0 & 100 &     4 &   214 &  78.8200 &     7 &    6495 & 0.0522\% &        168 &      3 &      242 \\
\bottomrule
\end{tabular}

    \caption{Dataset \textit{Glass}, $k$=6.}
    \label{tab:glass_detailed}
\end{table}

\begin{table}[!h]
\centering
\footnotesize
    \begin{tabular}{ccccccccccc}
\toprule
 ml &  cl &  seed &  size &       f &  cp\textsubscript{0} &  ineq\textsubscript{0} &   gap\textsubscript{0} & time\textsubscript{0}(s) &  nodes &  time(s) \\
\midrule
100 &   0 &     0 &   199 & 3336.54 &     8 &    6448 & 0.0002\% &        181 &      1 &      184 \\
100 &   0 &     1 &   199 & 3352.49 &     8 &    5856 & 0.0000\% &        377 &      1 &      380 \\
100 &   0 &     2 &   199 & 3337.51 &     9 &    3019 & 0.0002\% &        349 &      1 &      352 \\
100 &   0 &     3 &   199 & 3331.70 &     8 &    6019 & 0.0003\% &        390 &      1 &      393 \\
100 &   0 &     4 &   199 & 3319.86 &     9 &    6842 & 0.0002\% &        234 &      1 &      237 \\\hline
 50 &  50 &     0 &   236 & 3333.49 &     7 &    6901 & 0.0001\% &        197 &      1 &      199 \\
 50 &  50 &     1 &   231 & 3349.22 &    10 &    8046 & 0.0002\% &        336 &      1 &      339 \\
 50 &  50 &     2 &   225 & 3332.96 &     7 &    6520 & 0.1026\% &        167 &      5 &      785 \\
 50 &  50 &     3 &   232 & 3328.87 &     9 &    6737 & 0.0002\% &        477 &      1 &      480 \\
 50 &  50 &     4 &   232 & 3342.07 &     8 &    8136 & 0.0002\% &        333 &      1 &      336 \\\hline
  0 & 100 &     0 &   241 & 3341.28 &    12 &    6609 & 0.0000\% &        717 &      1 &      720 \\
  0 & 100 &     1 &   240 & 3338.45 &     7 &    7820 & 0.0004\% &        348 &      1 &      351 \\
  0 & 100 &     2 &   238 & 3338.48 &    10 &    4485 & 0.0002\% &        756 &      1 &      759 \\
  0 & 100 &     3 &   238 & 3354.66 &     8 &    7571 & 0.0002\% &        364 &      1 &      367 \\
  0 & 100 &     4 &   241 & 3322.84 &     7 &    7114 & 0.0002\% &        209 &      1 &      212 \\\hline
150 &   0 &     0 &   153 & 3364.21 &     5 &    6084 & 0.0000\% &        162 &      1 &      164 \\
150 &   0 &     1 &   152 & 3347.59 &     5 &    3358 & 0.0000\% &        180 &      1 &      182 \\
150 &   0 &     2 &   154 & 3335.71 &     6 &    5383 & 0.0003\% &        195 &      1 &      198 \\
150 &   0 &     3 &   155 & 3356.64 &     4 &    5802 & 0.0000\% &         70 &      1 &       72 \\
150 &   0 &     4 &   150 & 3362.72 &     5 &    6103 & 0.0001\% &        116 &      1 &      119 \\\hline
 75 &  75 &     0 &   188 & 3344.25 &     4 &    5664 & 0.0029\% &        127 &      1 &      130 \\
 75 &  75 &     1 &   191 & 3353.23 &     6 &    6434 & 0.0000\% &        186 &      1 &      188 \\
 75 &  75 &     2 &   187 & 3366.68 &     5 &    5425 & 0.0002\% &        313 &      1 &      316 \\
 75 &  75 &     3 &   194 & 3356.26 &     5 &    6524 & 0.0001\% &        295 &      1 &      298 \\
 75 &  75 &     4 &   187 & 3361.72 &     4 &    5998 & 0.0002\% &        155 &      1 &      158 \\\hline
  0 & 150 &     0 &   173 & 3369.50 &     5 &    2457 & 0.0002\% &        323 &      1 &      326 \\
  0 & 150 &     1 &   190 & 3357.44 &     7 &    3630 & 0.0002\% &        255 &      1 &      258 \\
  0 & 150 &     2 &   183 & 3365.68 &     5 &    3162 & 0.0001\% &        269 &      1 &      272 \\
  0 & 150 &     3 &   182 & 3356.13 &     5 &    5729 & 0.0010\% &        196 &      1 &      199 \\
  0 & 150 &     4 &   184 & 3357.68 &     5 &    6063 & 0.0000\% &        106 &      1 &      108 \\
\bottomrule
\end{tabular}

    \caption{Dataset \textit{Heart}, $k$=2.}
    \label{tab:heart_detailed}
\end{table}

\begin{table}[!h]
\centering
\footnotesize
    \begin{tabular}{ccccccccccc}
\toprule
 ml &  cl &  seed &  size &       f &  cp\textsubscript{0} &  ineq\textsubscript{0} &   gap\textsubscript{0} & time\textsubscript{0}(s) &  nodes &  time(s) \\
\midrule
100 &   0 &     0 &   212 & 552058 &    15 &    9326 & 0.0174\% &        744 &      3 &      974 \\
100 &   0 &     1 &   210 & 563454 &     9 &   10098 & 0.0229\% &        345 &      9 &      636 \\
100 &   0 &     2 &   211 & 538309 &     8 &   11644 & 0.0189\% &        308 &      5 &      461 \\
100 &   0 &     3 &   210 & 567156 &     5 &    8848 & 0.0096\% &        411 &      1 &      413 \\
100 &   0 &     4 &   210 & 533412 &     8 &   14579 & 0.0215\% &        295 &      7 &      490 \\\hline
 50 &  50 &     0 &   242 & 540393 &    12 &    6696 & 0.0006\% &        404 &      1 &      407 \\
 50 &  50 &     1 &   240 & 524886 &     5 &    7678 & 0.0011\% &        152 &      1 &      155 \\
 50 &  50 &     2 &   249 & 531530 &    10 &    8682 & 0.0048\% &        377 &      1 &      379 \\
 50 &  50 &     3 &   245 & 539132 &     7 &    8610 & 0.0327\% &        263 &     13 &      697 \\
 50 &  50 &     4 &   244 & 547714 &     7 &    6281 & 0.0106\% &        287 &      1 &      290 \\\hline
  0 & 100 &     0 &   252 & 548778 &     9 &    7476 & 0.0087\% &        432 &      1 &      434 \\
  0 & 100 &     1 &   257 & 544749 &     9 &    7391 & 0.0049\% &        324 &      1 &      327 \\
  0 & 100 &     2 &   252 & 538657 &     8 &    8410 & 0.0100\% &        322 &      1 &      325 \\
  0 & 100 &     3 &   256 & 532571 &     9 &    7834 & 0.0020\% &        488 &      1 &      491 \\
  0 & 100 &     4 &   257 & 528318 &     9 &    4842 & 0.0131\% &        460 &      3 &      548 \\\hline
150 &   0 &     0 &   161 & 605694 &     8 &    8834 & 0.1775\% &        325 &     13 &      609 \\
150 &   0 &     1 &   162 & 584261 &     2 &    7638 & 0.0120\% &         36 &      3 &       68 \\
150 &   0 &     2 &   161 & 575488 &     4 &    7159 & 0.0149\% &         63 &      3 &       87 \\
150 &   0 &     3 &   163 & 596006 &     4 &    8288 & 0.0080\% &         61 &      1 &       63 \\
150 &   0 &     4 &   160 & 592318 &     6 &    9186 & 0.0395\% &        382 &      3 &      477 \\\hline
 75 &  75 &     0 &   205 & 580718 &     4 &   10844 & 0.0142\% &        140 &      3 &      174 \\
 75 &  75 &     1 &   198 & 583640 &     7 &    8379 & 0.0204\% &        234 &      3 &      274 \\
 75 &  75 &     2 &   199 & 577566 &     5 &    7823 & 0.0200\% &        143 &      7 &      229 \\
 75 &  75 &     3 &   196 & 594787 &     2 &    7880 & 0.0066\% &         56 &      1 &       58 \\
 75 &  75 &     4 &   198 & 584670 &     4 &    9096 & 0.0214\% &        171 &      5 &      269 \\\hline
  0 & 150 &     0 &   197 & 568118 &     5 &    9549 & 0.0233\% &        168 &      3 &      234 \\
  0 & 150 &     1 &   200 & 554651 &     5 &    7186 & 0.0149\% &        117 &      3 &      147 \\
  0 & 150 &     2 &   197 & 562492 &     7 &    7015 & 0.0606\% &        194 &     29 &      609 \\
  0 & 150 &     3 &   201 & 576285 &     5 &    9863 & 0.0222\% &        225 &     13 &      593 \\
  0 & 150 &     4 &   201 & 565172 &     5 &    7615 & 0.0300\% &        136 &     13 &      352 \\
\bottomrule
\end{tabular}

    \caption{Dataset \textit{Vertebral}, $k$=2.}
    \label{tab:vertebral_detailed}
\end{table}

\begin{table}[!h]
\centering
\footnotesize
    \begin{tabular}{ccccccccccc}
\toprule
 ml &  cl &  seed &  size &       f &  cp\textsubscript{0} &  ineq\textsubscript{0} &   gap\textsubscript{0} & time\textsubscript{0}(s) &  nodes &  time(s) \\
\midrule
100 &   0 &     0 &   230 & 35191.0 &     8 &    2985 & 0.0164\% &        141 &      3 &      249 \\
100 &   0 &     1 &   229 & 36093.7 &    11 &   11928 & 0.2191\% &        297 &      9 &      766 \\
100 &   0 &     2 &   230 & 37363.4 &     7 &    7224 & 0.0028\% &        424 &      1 &      427 \\
100 &   0 &     3 &   229 & 37940.2 &     9 &    5590 & 0.0015\% &        246 &      1 &      249 \\
100 &   0 &     4 &   229 & 35165.5 &     8 &    2928 & 0.0039\% &        161 &      1 &      164 \\\hline
 50 &  50 &     0 &   280 & 28695.2 &    12 &    5365 & 0.0197\% &        298 &     15 &      715 \\
 50 &  50 &     1 &   279 & 29047.5 &    10 &    8920 & 0.0256\% &        324 &      3 &      518 \\
 50 &  50 &     2 &   279 & 29424.3 &    10 &    4748 & 0.0096\% &        298 &      1 &      301 \\
 50 &  50 &     3 &   279 & 30652.2 &    13 &    7428 & 0.1675\% &        194 &     45 &     1622 \\
 50 &  50 &     4 &   279 & 26835.9 &     8 &    8879 & 0.0269\% &        280 &     17 &      761 \\\hline
  0 & 100 &     0 &   329 & 18796.9 &    10 &    5106 & 0.0117\% &        327 &      5 &      575 \\
  0 & 100 &     1 &   329 & 19050.2 &    13 &    7343 & 0.1368\% &        316 &     51 &     1259 \\
  0 & 100 &     2 &   329 & 18803.8 &    11 &    6707 & 0.0203\% &        339 &      7 &      800 \\
  0 & 100 &     3 &   329 & 18965.5 &    11 &    7663 & 0.1467\% &        274 &    200 (0.013\%) &     3008 \\
  0 & 100 &     4 &   329 & 18808.6 &    12 &    4871 & 0.0204\% &        343 &      7 &      680 \\\hline
150 &   0 &     0 &   184 & 42917.2 &     5 &    6023 & 0.8044\% &        617 &      3 &      824 \\
150 &   0 &     1 &   192 & 41496.5 &     6 &    6923 & 0.0020\% &        269 &      1 &      272 \\
150 &   0 &     2 &   182 & 42197.4 &     8 &    7770 & 1.0198\% &        147 &      9 &      414 \\
150 &   0 &     3 &   183 & 41485.8 &     7 &    7047 & 0.0182\% &        736 &      5 &     1124 \\
150 &   0 &     4 &   189 & 40234.0 &     8 &    7108 & 0.0234\% &        402 &      7 &      764 \\\hline
 75 &  75 &     0 &   254 & 29027.3 &    11 &    7667 & 0.2384\% &        164 &    200 (0.027\%) &     1308 \\
 75 &  75 &     1 &   254 & 30703.0 &    11 &    7529 & 0.1313\% &        188 &     51 &      740 \\
 75 &  75 &     2 &   254 & 33039.1 &    10 &    9060 & 0.0278\% &        230 &     33 &      788 \\
 75 &  75 &     3 &   255 & 35238.8 &     9 &   10840 & 1.8065\% &        115 &    200 (0.032\%) &     1386 \\
 75 &  75 &     4 &   255 & 33575.1 &     7 &   12486 & 0.0070\% &        167 &      1 &      170 \\\hline
  0 & 150 &     0 &   329 & 19194.0 &    15 &    8291 & 0.2900\% &        281 &    125 &     1817 \\
  0 & 150 &     1 &   329 & 19127.6 &    12 &    7877 & 0.0739\% &        267 &     19 &     1017 \\
  0 & 150 &     2 &   329 & 19117.6 &    14 &    7485 & 0.0941\% &        299 &    207 &     2686 \\
  0 & 150 &     3 &   329 & 18995.9 &    12 &    5656 & 0.0239\% &        321 &     15 &      698 \\
  0 & 150 &     4 &   329 & 19181.9 &    13 &    7281 & 0.0540\% &        321 &      9 &      747 \\
\bottomrule
\end{tabular}

    \caption{Dataset \textit{Accent}, $k$=6.}
    \label{tab:accent_detailed}
\end{table}

\begin{table}[!h]
\centering
\footnotesize
    \begin{tabular}{ccccccccccc}
\toprule
 ml &  cl &  seed &  size &       f &  cp\textsubscript{0} &  ineq\textsubscript{0} &   gap\textsubscript{0} & time\textsubscript{0}(s) &  nodes &  time(s) \\
\midrule
100 &   0 &     0 &   236 & 16.6586 &     6 &    8249 & 0.0137\% &        146 &      5 &      249 \\
100 &   0 &     1 &   236 & 15.8052 &     5 &    6358 & 0.0259\% &        155 &     43 &     1117 \\
100 &   0 &     2 &   237 & 17.2182 &     5 &    5993 & 0.0799\% &        151 &     59 &     1010 \\
100 &   0 &     3 &   236 & 17.3390 &     9 &    5781 & 0.1500\% &        254 &     35 &     1192 \\
100 &   0 &     4 &   237 & 16.7166 &     7 &    6787 & 0.0277\% &        187 &     25 &      692 \\\hline
 50 &  50 &     0 &   286 & 14.8403 &     6 &    6485 & 0.0978\% &        241 &     11 &      633 \\
 50 &  50 &     1 &   286 & 14.9577 &     8 &    7541 & 0.0412\% &        237 &     17 &      755 \\
 50 &  50 &     2 &   286 & 15.9729 &     6 &    7198 & 0.0111\% &        188 &     19 &      528 \\
 50 &  50 &     3 &   286 & 15.1497 &     9 &    8834 & 0.0602\% &        265 &      3 &      426 \\
 50 &  50 &     4 &   286 & 15.8739 &     8 &    7708 & 0.0381\% &        186 &     85 &     1462 \\\hline
  0 & 100 &     0 &   336 & 13.8480 &    11 &    8135 & 0.0640\% &        331 &    129 &     2591 \\
  0 & 100 &     1 &   336 & 13.9326 &     8 &    7998 & 0.0887\% &        256 &    200 (0.024\%) &     3054 \\
  0 & 100 &     2 &   336 & 13.9852 &     9 &    8319 & 0.0435\% &        294 &    109 &     2029 \\
  0 & 100 &     3 &   336 & 14.0955 &    11 &    9490 & 0.2282\% &        291 &    200 (0.023\%) &     2962 \\
  0 & 100 &     4 &   336 & 13.9024 &     8 &    6923 & 0.0971\% &        266 &     83 &     1524 \\\hline
150 &   0 &     0 &   189 & 16.6943 &     7 &   13495 & 0.0112\% &        693 &      5 &     1184 \\
150 &   0 &     1 &   186 & 17.0550 &     4 &    6228 & 0.0074\% &        162 &      1 &      165 \\
150 &   0 &     2 &   189 & 17.6560 &     6 &    8364 & 0.0083\% &        543 &      1 &      546 \\
150 &   0 &     3 &   194 & 16.9603 &     5 &    7442 & 0.0343\% &        315 &      5 &      552 \\
150 &   0 &     4 &   189 & 17.7616 &     5 &    6345 & 0.0141\% &        158 &      5 &      344 \\\hline
 75 &  75 &     0 &   262 & 16.1113 &     6 &    7573 & 0.0166\% &        197 &     17 &      649 \\
 75 &  75 &     1 &   261 & 17.0244 &     7 &    7595 & 0.7033\% &        149 &     91 &     1051 \\
 75 &  75 &     2 &   262 & 16.3000 &     8 &    6991 & 0.6547\% &        189 &     77 &     1052 \\
 75 &  75 &     3 &   261 & 16.0741 &     5 &    8682 & 0.0416\% &        701 &      3 &      779 \\
 75 &  75 &     4 &   261 & 16.0650 &     8 &    7768 & 0.0594\% &        213 &     49 &     1030 \\\hline
  0 & 150 &     0 &   336 & 14.1895 &     8 &    8102 & 0.4973\% &        201 &     83 &     1844 \\
  0 & 150 &     1 &   336 & 14.1301 &    14 &    5932 & 0.1828\% &        313 &     29 &     1008 \\
  0 & 150 &     2 &   336 & 14.2708 &    10 &    9263 & 0.1760\% &        348 &    199 &     3432 \\
  0 & 150 &     3 &   336 & 14.0998 &    11 &   10039 & 0.1347\% &        301 &    113 &     2263 \\
  0 & 150 &     4 &   336 & 13.9936 &     8 &    7811 & 0.3591\% &        226 &    159 &     2612 \\
\bottomrule
\end{tabular}

    \caption{Dataset \textit{Ecoli}, $k$=8.}
    \label{tab:ecoli_detailed}
\end{table}

\begin{table}[!h]
\centering
\footnotesize
    \begin{tabular}{ccccccccccc}
\toprule
 ml &  cl &  seed &  size &       f &  cp\textsubscript{0} &  ineq\textsubscript{0} &   gap\textsubscript{0} & time\textsubscript{0}(s) &  nodes &  time(s) \\
\midrule
150 &   0 &     0 &   350 & 13157.1 &    10 &    9700 & 0.0105\% &        655 &      3 &      842 \\
150 &   0 &     1 &   350 & 13063.7 &    16 &    8438 & 0.0068\% &        895 &      1 &      899 \\
150 &   0 &     2 &   351 & 13141.3 &    11 &    7749 & 0.0230\% &        448 &     23 &     1103 \\
150 &   0 &     3 &   350 & 13504.3 &    18 &    8176 & 0.0125\% &        771 &     11 &     1519 \\
150 &   0 &     4 &   351 & 13302.7 &    11 &    7926 & 0.0154\% &        548 &      7 &     1069 \\\hline
 75 &  75 &     0 &   426 & 12644.9 &    19 &    9083 & 0.0377\% &       1052 &     15 &     2561 \\
 75 &  75 &     1 &   425 & 12684.0 &    22 &    9297 & 0.0457\% &       1133 &    123 &     5470 \\
 75 &  75 &     2 &   425 & 12880.3 &    21 &    9778 & 0.0732\% &       1327 &      3 &     2553 \\
 75 &  75 &     3 &   425 & 12397.3 &    18 &    5537 & 0.0384\% &       1008 &     11 &     2259 \\
 75 &  75 &     4 &   425 & 12750.0 &    21 &    8853 & 0.0238\% &       1417 &     17 &     3850 \\\hline
  0 & 150 &     0 &   500 & 11615.7 &    23 &    5724 & 0.0089\% &       2181 &      1 &     2186 \\
  0 & 150 &     1 &   500 & 11681.6 &    23 &   10838 & 0.0624\% &       1796 &     75 &     7566 \\
  0 & 150 &     2 &   500 & 11737.5 &    26 &   10256 & 0.0271\% &       2421 &      5 &     4555 \\
  0 & 150 &     3 &   500 & 11690.5 &    28 &    8392 & 0.0275\% &       2452 &    200 (0.012\%) &    15416 \\
  0 & 150 &     4 &   500 & 11613.6 &    23 &    6071 & 0.0332\% &       2066 &     35 &    10491 \\\hline
250 &   0 &     0 &   252 & 14234.3 &     8 &    9178 & 0.0501\% &        246 &      5 &      427 \\
250 &   0 &     1 &   253 & 14096.2 &    10 &    7906 & 1.4892\% &        483 &      7 &      922 \\
250 &   0 &     2 &   254 & 14113.0 &     6 &    6366 & 0.0317\% &        267 &     11 &      772 \\
250 &   0 &     3 &   255 & 13594.2 &     6 &    8954 & 0.0104\% &        386 &      3 &      542 \\
250 &   0 &     4 &   252 & 14297.6 &     8 &    7014 & 0.0387\% &        424 &      5 &      733 \\\hline
125 & 125 &     0 &   375 & 13127.4 &    15 &    3903 & 0.0152\% &        812 &      5 &     1110 \\
125 & 125 &     1 &   375 & 13129.9 &    20 &    5653 & 0.0728\% &        888 &      3 &     1454 \\
125 & 125 &     2 &   375 & 12656.5 &    12 &    7445 & 0.0031\% &       1097 &      1 &     1101 \\
125 & 125 &     3 &   375 & 13034.5 &    19 &    5570 & 0.2095\% &        882 &     39 &     2590 \\
125 & 125 &     4 &   375 & 13049.9 &    16 &    4840 & 0.0437\% &        824 &     97 &     4184 \\\hline
  0 & 250 &     0 &   500 & 11791.9 &    23 &   10565 & 0.0957\% &       1791 &     45 &     8007 \\
  0 & 250 &     1 &   500 & 11760.4 &    30 &    6126 & 0.0025\% &       2811 &      1 &     2817 \\
  0 & 250 &     2 &   500 & 11685.1 &    19 &    9463 & 0.0318\% &       1468 &      7 &     3136 \\
  0 & 250 &     3 &   500 & 11758.2 &    24 &    8575 & 0.0036\% &       2322 &      1 &     2327 \\
  0 & 250 &     4 &   500 & 11846.7 &    30 &   10099 & 0.0266\% &       2775 &     59 &     9984 \\
\bottomrule
\end{tabular}

    \caption{Dataset \textit{ECG5000}, $k$=5.}
    \label{tab:ECG5000_detailed}
\end{table}

\begin{table}[!h]
\centering
\footnotesize
    \begin{tabular}{ccccccccccc}
\toprule
 ml &  cl &  seed &  size &       f &  cp\textsubscript{0} &  ineq\textsubscript{0} &   gap\textsubscript{0} & time\textsubscript{0}(s) &  nodes &  time(s) \\
\midrule
150 &   0 &     0 &   350 & 315572 &    13 &    7511 & 0.0002\% &       1376 &      1 &     1380 \\
150 &   0 &     1 &   350 & 316777 &    16 &    4452 & 0.0003\% &       1479 &      1 &     1483 \\
150 &   0 &     2 &   351 & 317242 &    18 &    7530 & 0.0004\% &       1707 &      1 &     1710 \\
150 &   0 &     3 &   351 & 316870 &    16 &    7341 & 0.0135\% &        998 &      3 &     1402 \\
150 &   0 &     4 &   350 & 319417 &    25 &    5954 & 0.0004\% &       2283 &      1 &     2286 \\\hline
 75 &  75 &     0 &   411 & 317837 &    19 &    6894 & 0.0005\% &       1571 &      1 &     1575 \\
 75 &  75 &     1 &   402 & 317099 &    16 &    6948 & 0.0002\% &       2546 &      1 &     2550 \\
 75 &  75 &     2 &   405 & 317551 &    16 &    7938 & 0.0004\% &       1790 &      1 &     1793 \\
 75 &  75 &     3 &   404 & 316084 &    12 &    7576 & 0.0001\% &       1338 &      1 &     1341 \\
 75 &  75 &     4 &   398 & 318000 &    25 &    8418 & 0.0594\% &       2350 &      3 &     2967 \\\hline
  0 & 150 &     0 &   425 & 317954 &    18 &    7900 & 0.0002\% &       1869 &      1 &     1873 \\
  0 & 150 &     1 &   425 & 317061 &    19 &    6704 & 0.0003\% &       2084 &      1 &     2088 \\
  0 & 150 &     2 &   432 & 317114 &    20 &    8044 & 0.0001\% &       2797 &      1 &     2801 \\
  0 & 150 &     3 &   427 & 317781 &    20 &    8547 & 0.0021\% &       2529 &      1 &     2533 \\
  0 & 150 &     4 &   427 & 318294 &    19 &    8152 & 0.0002\% &       2457 &      1 &     2461 \\\hline
250 &   0 &     0 &   251 & 323167 &     6 &    6592 & 0.0003\% &        187 &      1 &      190 \\
250 &   0 &     1 &   253 & 324203 &    10 &    2839 & 0.0004\% &        774 &      1 &      777 \\
250 &   0 &     2 &   250 & 322362 &     8 &    2069 & 0.0007\% &        247 &      1 &      250 \\
250 &   0 &     3 &   253 & 323319 &    10 &    6509 & 0.0006\% &        384 &      1 &      387 \\
250 &   0 &     4 &   253 & 322140 &     7 &    2650 & 0.0004\% &        332 &      1 &      334 \\\hline
125 & 125 &     0 &   324 & 325922 &    14 &    4988 & 0.0004\% &       1585 &      1 &     1588 \\
125 & 125 &     1 &   327 & 323445 &    10 &    6743 & 0.0005\% &       1075 &      1 &     1078 \\
125 & 125 &     2 &   318 & 326321 &    14 &    8061 & 0.0005\% &       2111 &      1 &     2114 \\
125 & 125 &     3 &   325 & 323682 &     8 &    6815 & 0.0001\% &        360 &      1 &      363 \\
125 & 125 &     4 &   322 & 322687 &    10 &    2413 & 0.0002\% &        684 &      1 &      687 \\\hline
  0 & 250 &     0 &   318 & 324062 &    10 &    6619 & 0.0002\% &        733 &      1 &      736 \\
  0 & 250 &     1 &   329 & 323982 &     8 &    7062 & 0.0001\% &        615 &      1 &      618 \\
  0 & 250 &     2 &   320 & 324604 &     8 &    7116 & 0.0003\% &        514 &      1 &      517 \\
  0 & 250 &     3 &   322 & 325706 &    11 &    7502 & 0.0008\% &        493 &      1 &      496 \\
  0 & 250 &     4 &   317 & 325240 &    15 &    7276 & 0.0007\% &        842 &      1 &      845 \\
\bottomrule
\end{tabular}

    \caption{Dataset \textit{Computers}, $k$=2.}
    \label{tab:computers_detailed}
\end{table}

\begin{table}[!h]
\centering
\footnotesize
    \begin{tabular}{ccccccccccc}
\toprule
 ml &  cl &  seed &  size &       f &  cp\textsubscript{0} &  ineq\textsubscript{0} &   gap\textsubscript{0} & time\textsubscript{0}(s) &  nodes &  time(s) \\
\midrule
200 &   0 &     0 &   604 & 1.78019e+07 &     2 &    6361 & 0.0013\% &        505 &        1 &      512 \\
200 &   0 &     1 &   601 & 1.78183e+07 &     1 &    5639 & 0.0002\% &        439 &        1 &      446 \\
200 &   0 &     2 &   601 & 1.78077e+07 &     2 &    6629 & 0.0019\% &        516 &        1 &      522 \\
200 &   0 &     3 &   601 & 1.78019e+07 &     2 &    6987 & 0.0019\% &        524 &        1 &      531 \\
200 &   0 &     4 &   601 & 1.78082e+07 &     2 &    8324 & 0.0004\% &        501 &        1 &      508 \\\hline
100 & 100 &     0 &   701 & 1.78171e+07 &     2 &    7395 & 0.0059\% &        660 &        1 &      670 \\
100 & 100 &     1 &   701 & 1.78019e+07 &     2 &    8785 & 0.0057\% &        690 &        1 &      700 \\
100 & 100 &     2 &   701 & 1.78132e+07 &     2 &    7318 & 0.0026\% &        623 &        1 &      633 \\
100 & 100 &     3 &   701 & 1.78257e+07 &     2 &    8034 & 0.0153\% &        718 &        7 &     2608 \\
100 & 100 &     4 &   701 & 1.78064e+07 &     2 &    7866 & 0.0117\% &        669 &        3 &     1317 \\\hline
  0 & 200 &     0 &   801 & 1.78022e+07 &     2 &   10000 & 0.0094\% &        859 &        1 &      873 \\
  0 & 200 &     1 &   801 & 1.78019e+07 &     3 &    9779 & 0.0119\% &        918 &        3 &     1764 \\
  0 & 200 &     2 &   801 & 1.78031e+07 &     2 &   10000 & 0.0133\% &        828 &        5 &     2407 \\
  0 & 200 &     3 &   801 & 1.78027e+07 &     2 &   10000 & 0.0159\% &        758 &        3 &     1567 \\
  0 & 200 &     4 &   801 & 1.78019e+07 &     5 &   11572 & 0.0068\% &       1075 &        1 &     1089 \\\hline
400 &   0 &     0 &   413 & 1.78171e+07 &     2 &   10402 & 0.0002\% &        345 &        1 &      349 \\
400 &   0 &     1 &   419 & 1.78067e+07 &     2 &    5317 & 0.0001\% &        649 &        1 &      653 \\
400 &   0 &     2 &   419 & 1.78130e+07 &     1 &    5837 & 0.0000\% &        337 &        1 &      341 \\
400 &   0 &     3 &   414 & 1.78073e+07 &     1 &    6059 & 0.0003\% &        311 &        1 &      315 \\
400 &   0 &     4 &   421 & 1.78125e+07 &     3 &    3448 & 0.0001\% &       1091 &        1 &     1095 \\\hline
200 & 200 &     0 &   602 & 1.78137e+07 &     2 &    6452 & 0.0208\% &        537 &        9 &     2638 \\
200 & 200 &     1 &   603 & 1.78023e+07 &     2 &    6083 & 0.0002\% &        446 &        1 &      453 \\
200 & 200 &     2 &   602 & 1.78026e+07 &     2 &    5327 & 0.0002\% &        519 &        1 &      526 \\
200 & 200 &     3 &   602 & 1.78125e+07 &     3 &    9862 & 0.0009\% &        527 &        1 &      534 \\
200 & 200 &     4 &   602 & 1.78168e+07 &     4 &    8226 & 0.0074\% &        629 &        1 &      636 \\\hline
  0 & 400 &     0 &   801 & 1.78064e+07 &     4 &   10125 & 0.0031\% &       1111 &        1 &     1124 \\
  0 & 400 &     1 &   801 & 1.78064e+07 &     3 &    8839 & 0.0171\% &        907 &        3 &     1809 \\
  0 & 400 &     2 &   801 & 1.78131e+07 &     3 &    9704 & 0.0179\% &        949 &        3 &     1880 \\
  0 & 400 &     3 &   801 & 1.78068e+07 &     2 &    9638 & 0.0174\% &        784 &        5 &     2434 \\
  0 & 400 &     4 &   801 & 1.78026e+07 &     2 &   10000 & 0.0120\% &        760 &        3 &     1593 \\
\bottomrule
\end{tabular}

    \caption{Dataset \textit{Gene}, $k$=5.}
    \label{tab:gene_detailed}
\end{table}

\clearpage

\bibliography{sample}

\end{document}